\documentclass[microtype]{gtpart}     
\agtart

\usepackage{amscd}

\title{Exponential growth of torsion in abelian coverings}

\author{Jean Raimbault}
\givenname{Jean}
\surname{Raimbault}
\address{Institut de math\'ematiques de Jussieu\\ 4, place Jussieu \\Case 247\\ 75252 Paris C\'edex 5\\ France}
\email{raimbault@math.jussieu.fr}
\urladdr{http://www.math.jussieu.fr/~raimbault/}

\keyword{Reidemeister torsion}
\keyword{$\ell^2$-torsion}
\subject{primary}{msc2000}{57M10}
\subject{secondary}{msc2000}{57M25}
\subject{secondary}{msc2000}{57Q10}

\arxivreference{math.GT/1012.3666}
\arxivpassword{qex89}

\volumenumber{}
\issuenumber{}
\publicationyear{}
\papernumber{}
\startpage{}
\endpage{}
\doi{}
\MR{}
\Zbl{}
\received{}
\revised{}
\accepted{}
\published{}
\publishedonline{}
\proposed{}
\seconded{}
\corresponding{}
\editor{}
\version{}

\newtheorem{theo}{Theorem}[section]    
\newtheorem{lem}[theo]{Lemma}          
\newtheorem{prop}[theo]{Proposition}
\newtheorem*{cor}{Corollary}

\newcommand{\vol}{\operatorname{vol}}
\newcommand{\im}{\operatorname{Im}}

\newcommand{\rg}{\operatorname{rk}}
\newcommand{\coker}{\operatorname{coker}}
\newcommand{\pdet}{\det{}'}
\newcommand{\fkdet}{\det{}_{\mathcal{N}(G)}}

\begin{document}

\begin{abstract}    
We show exponential growth of torsion numbers for links whose first
nonzero Alexander polynomial has positive logarithmic Mahler measure. 
This extends a theorem of Silver and Williams to the case
of a null first Alexander polynomial and provides a partial solution
for a conjecture of theirs.
\end{abstract}

\maketitle


\numberwithin{equation}{section}

\section*{Introduction}

Let $M$ be a compact three--manifold; the homology groups $H_i(M)$ can be 
written as the direct sums $H_i(M)_{tors}\oplus H_i(M)_{free}$ of a finite 
abelian group with a finite-rank free abelian group. The torsion summand is 
nontrivial only for $i=1$: $H_0$ and $H_3$ are $\mathbb{Z}$ or $0$, and the 
universal coefficients theorem for cohomology implies that $H^1(M,\partial M)$ 
is free, and by Poincar\'e duality it follows that $H_2(M)$ is also 
torsion-free. On the other hand, the torsion in $H_1$ can be arbitrarily large
(e.g., for lens spaces; see below for hyperbolic examples) and it is believed 
that ``most'' 3-manifolds should have a rather large torsion. For example 
E. Kowalski shows in \cite[Proposition 7.19]{Kowalski} that the first homology 
group of a ``Dunfield-Thurston random 3-manifold'' typically has a large 
torsion subgroup. This paper is concerned with the growth rate of the order of 
$H_1(M_N)_{tors}$ in a sequence of finite coverings $M_N$ of a
manifold $M$. The least precise question that can be asked is whether
it is exponential in the degree or not, i.e. whether the sequence 
$\frac{\log|H_1(M_N)_{tors}|}{[\pi_1(M):\pi_1(M_N)]}$ has a positive limit 
(or limit superior). This shall be partially answered here in the case where 
the $M_N$ are abelian coverings converging to a free abelian covering of $M$. 
The main motivation to study this question was to provide a partial result 
towards a conjecture of Silver and Williams on the growth rate of torsion 
numbers of abelian coverings of complements of links (see Conjecture 6.1 in 
\cite{Silver/Williams1} or \eqref{conjecture} below).

Historically, the first context where the 
growth of torsion in the homology of coverings has been studied is that of 
cyclic coverings of a knot complement. Let $K$ be an open knotted solid torus 
in the three-sphere and $M=\mathbb{S}^3-K$; then $M$ is a compact 
three-manifold with $H_1(M)=\mathbb{Z}$. Thus we can consider the infinite 
cyclic covering $\widehat{M}$ of $M$ and its finite quotients 
$M_N$, which are the finite coverings of $M$ corresponding to the maps 
$\pi_1(M)\rightarrow \mathbb{Z}\rightarrow \mathbb{Z}/N\mathbb{Z}$. If 
$\Delta$ is the Alexander polynomial of $K$ (the first Alexander polynomial 
of the $\mathbb{Z}[\mathbb{Z}]$-module $H_1(\widehat{M})$-see 
\ref{Alexander1}) then 
the so-called Fox formula says that for all $N$ such that $H_1(M_N)$ has rank 
one we have:
$$|H_1(M_N)_{tors}|=\prod_{\zeta^N=1}|\Delta(\zeta)|.$$
It is known that $\Delta$ is always nonzero. A diophantine 
inequality due to Gelfond then allows to show that when $N$ tends 
to infinity  $\frac 1 N\sum_{\zeta^N=1}\log|\Delta(\zeta)|$ converges to 
the logarithmic Mahler measure of $\Delta$, that is:
$$m(\Delta)=\int_{\mathbb{T}^1}\log|\Delta(z)|dz$$
(here $\mathbb{T}^1$ is the unit circle in $\mathbb{C}$ endowed with the 
Lebesgue probability measure); see Lemma \ref{m=1} below or Proposition 2.8 of 
\cite{Silver/Williams2}. Let $d$ be the product of all integers $l$ such that 
the $l$th cyclotomic polynomial divides $\Delta$: it is known 
(\cite[Proposition 2.2]{Silver/Williams2} or \cite[Theorem 4.1(ii)]{Gordon}) 
that $H_1(M_N,\mathbb{Q})\cong H_1(M_{N+d},\mathbb{Q})$. We get that
$$\lim_{\substack{N\rightarrow\infty\\ d|N-1}}(\frac 1 N \log|H_1(M_N)_{tors}|)=m(\Delta);$$
in particular, if $m(\Delta)>0$ an infinite subsequence of the torsion numbers 
tends to infinity with exponential growth. This result 
answered a question of Gordon (who proved in \cite{Gordon} that when 
$m(\Delta)=0$ the torsion numbers are periodic in $N$) and has first been 
proved by Gonz{\'a}lez-Acu{\~n}a and Short in \cite{G-A/S} and independently 
by Robert Riley in \cite{Riley}. Since there exists hyperbolic knots whose 
Alexander polynomial has positive Mahler measure (for example the figure-eight)
 this proves in particular that there exists (noncompact with finite volume) 
hyperbolic manifolds with arbitrarily large torsion in their $H_1$.

The similar setting for a collection of $m>1$ linked solid tori $L$ in the 
three-sphere is to consider the coverings $M_H$ of $M=\mathbb{S}^3-L$ given 
by the maps $\pi_1(M)\rightarrow H_1(M)\cong\mathbb{Z}^m\rightarrow 
\mathbb{Z}^m/H$ where $H$ is a finite-index subgroup. We want to study the 
torsion numbers as the quantity:
$$\alpha(H)=\min\{\max_i|v_i|;\: 0\not=v=(v_1,\ldots,v_m)\in H\}$$
tends to infinity. Fox's formula admits a generalization due to Mayberry 
and Murasugi (\cite{Mayberry/Murasugi}; see also \cite{Porti} for a proof 
using Reidemeister torsion). However, because the Betti numbers may very well 
tend to infinity as $\alpha(H)$ does, it may not be applicable to an infinite 
sequence of coverings. This forbids any na\"ive application of the preceding 
scheme of proof to this case. Moreover, two 
other obstructions arise when considering links: the Alexander polynomial 
may be zero, and the generalisation of Gelfond's estimate needed to 
prove the convergence of Riemann sums to the Mahler measure is not known to 
hold for polynomials in several variables.

By using methods from algebraic dynamical systems Daniel Silver and Susan 
Williams, in \cite{Silver/Williams1}, were able to show that if the first 
Alexander polynomial $\Delta(L)$ of $L$ is nonzero then:
$$\limsup_{\alpha(H)\rightarrow\infty}\frac{\log|H_1(M_H)_{tors}|}{[G:H]}=m(\Delta(L)):=\int_{\mathbb{T}^m}\log|\Delta(L)|$$
and that this holds with a limit in the case of a knot (this is also proved 
in \cite{Silver/Williams2}). They also conjecture that, in the general case, 
the following limit should hold where $\Delta_i$ is the first nonzero 
polynomial in the sequence of Alexander polynomials of $L$:
\begin{equation}
\lim_{\alpha(H)\rightarrow\infty}\frac{\log|H_1(M_H)_{tors}|}{[G:H]}=m(\Delta_i(L)).
\label{conjecture}
\end{equation}

The Mahler measure of the Alexander polynomial of a knot can be interpreted as 
the $\ell^2$-torsion of the infinite cyclic covering $\widehat{M}$ of its 
exterior $M$ (see for example \eqref{L2tors} below). 
In \cite{Bergeron/Venkatesh} Nicolas Bergeron and Akshay Venkatesh used this 
to give a new proof of the theorem of Silver and Williams in the case of a 
knot: the main point is that the nonnullity of the Alexander polynomial 
guarantees that $\widehat{M}$ is $\ell^2$-acyclic, and in this case the 
approximation of $\ell^2$-torsion by Reidemeister torsions is known. They then 
show that the growth of torsion numbers is the same as that of 
Reidemeister torsion and obtain a generalization of the theorem for 
knots (see \cite[Theorem 7.3]{Bergeron/Venkatesh}). 
For links with nonzero first Alexander polynomial their 
proof can be adapted; however the result obtained is slightly weaker than 
Silver and Williams'.

The aim of this work is to generalise their results to the case of a null 
first Alexander polynomial, i.e. that of a non-$\ell^2$-acyclic infinite 
covering. This yields the following result, which follows from the more 
general Theorem \ref{approx}.

\begin{theo}
With notation as in \eqref{conjecture} we have:
$$\limsup_{\alpha(H)\rightarrow\infty}\frac{\log|H_1(M_H)_{tors}|}{[G:H]}\ge m(\Delta_i(L)).$$
\label{Main1}
\end{theo}

The  same scheme of proof can be applied to homology with coefficients in a 
local system. Let $V$ be a free 
$\mathbb{Z}$-module of finite rank and $\chi$ a representation of $\pi_1(M)$ 
to $GL(V)$. Then we can define homology twisted by $\chi$ for all coverings 
of $M$, also called homology with coefficients in the local system 
defined by $\chi$. See \cite[Chapter 3.H]{Hatcher} for definitions. 
Let $\Delta_i(L,\chi)$ be the first nonzero Alexander polynomial 
of $H_1(\widehat{M};V)$ (the Alexander polynomials with nontrivial 
coefficients are usually called twisted Alexander 
polynomials). We get (see also \cite[Theorem 3.10]{Silver/Williams3}):
$$\limsup_{\alpha(H)\rightarrow\infty}\frac{\log|H_1(M_H;V)_{tors}|}{[G:H]}\ge m(\Delta_i(L,\chi)).$$

In the case of a sequence of cyclic coverings converging to a possibly 
non-$\ell^2$-acyclic infinite cyclic covering the proofs yield a more precise 
result, generalising Theorem 2.10 of 
\cite{Silver/Williams2}. 

\begin{theo}
If $M$ is a compact $n$-manifold (with boundary) with a surjection 
$\pi_1(M)\rightarrow\mathbb{Z}$ let $M_N$ be the cyclic covering of $M$ 
associated to $\pi_1(M)\rightarrow\mathbb{Z}/N\mathbb{Z}$ and $\widehat{M}$ 
the infinite cyclic covering. Then we have for all $i=1,\ldots,n-1$:
$$\lim_{N\rightarrow\infty}\frac{\log|H_i(M_N)_{tors}|}{N}=m(\Delta_i)$$
where $\Delta_i$ is the first nonzero Alexander polynomial of 
$H_i(\widehat{M})$.
\label{Main2}
\end{theo}

Thang Le has recently proved the equality lacking in  Theorem \ref{Main1} 
in the preprint \cite{Le}. His proof uses methods slightly different 
from ours, namely a Bourbaki theorem on ``pseudo-zero modules''
and the dynamical version of approximation (due to Lind; see 
\cite[Theorem 21.1]{Schmidt}). He also gives a different construction of 
the approximating sequences contructed in Proposition \ref{explicit} below.

Finally, note that throughout this paper we work with three-manifolds 
with boundary, that is, we consider regular coverings of $\mathbb{S}^3-L$ 
instead of the correspondings branched (over core circles of $L$) coverings 
of $\mathbb{S}^3$, which are closed three-manifolds. However, the growth of 
torsion homology in abelian coverings is the same for the two sequences of 
manifolds, as shown in 
\cite[4.4]{Le} (for knots the two homologies differ 
only by a infinite cyclic direct factor).

\subsection*{Warnings and outline}
Most of the techniques and results used here are fairly elementary with one 
notable exception: I have found more practical and efficient to use the 
Fuglede-Kadison determinant as defined in the book \cite{Luck} of Wolfgang 
L\"uck, and all its properties proved therein instead of trying to do 
everything solely in terms of Mahler measures (in any case, it is the right 
tool to use to define $\ell^2$-torsion). However, the reader need not be 
familiar with the Fuglede-Kadison determinant since all those of its 
properties that we use are recalled. One warning about notations: I have 
chosen to denote by $\mathcal{M}$ the multiplicative Mahler measure given by 
$\exp(m)$ since I already use the casual $M$ to denote modules and manifolds. 

The paper is organised as follows: the first section reviews various 
classical and $\ell^2$-invariants for $\mathbb{Z}[G]$-modules. In particular, 
I give there the computation of the $\ell^2$-torsion in terms of Mahler 
measures of Alexander polynomials. In the second section I give a proof of 
the approximation of the Fuglede-Kadison determinant in the free Abelian 
case since there was (to the best of my knowledge) no complete proof 
of this result available in the literature. In the third section we prove 
Theorems \ref{Main1} and \ref{Main2}. The fourth section contains 
some further remarks about three-manifolds. The appendices contain proofs 
that I judged too cumbersome for a first reading; I hope that their being 
relegated there will make the main thrust of the paper more apparent and 
intelligible.

\subsection*{Acknowledgements} The subject treated here has been proposed to 
me by my Ph.D. advisor, Nicolas Bergeron. I gratefully thank him for this as 
well as for the numerous suggestions he made regarding the redaction of this 
paper and the math inside. I also wish to thank Pierre Vogel who tried to 
explain to me the Alexander polynomial. I have benefited from 
conversations with Thang Le regarding his paper \cite{Le} and related topics, 
as well as his pointing out a deadly mistake in a previous proof of 
Theorem \ref{torsion3}. Also, the Bombieri-Zannier theorem used in 
Appendix \ref{loose} has been brought to my attention by the notes from a 
lecture he gave. Finally, I thank Takahiro Kitayama for reading carefully a 
previous version of this paper and the referee for spotting many typos and 
minor flaws.


\section{Review of classical and $\ell^2$-invariants}

In this section $G$ shall denote a free Abelian group of rank $m$; we 
identify the group ring $\mathbb{Z}[G]$ with the Laurent polynomials 
$\mathbb{Z}[t_1^{\pm1},\ldots,t_m^{\pm1}]$ through the choice of a 
free generating set $t_1,\ldots,t_m$ for $G$. In this section, we review 
various invariants for $G$-modules and $G$-complexes.


\subsection{Rank, von Neumann dimension and $\ell^2$-Betti numbers}

Recall that a finitely generated Hilbert $\mathcal{N}(G)$-module is a Hilbert 
space with a unitary action of $G$, which admits a $G$-equivariant surjection 
of $\ell^2(G)^n$ for some integer $n$ (we shall not use explicitly the 
von Neumann algebra $\mathcal{N}(G)$ in the sequel, so we shall simply call 
these finitely generated $\mathcal{N}(G)$-modules). Given a submodule 
$M\subset\mathbb{Z}[G]^n$, the closure of $M\otimes\mathbb{C}$ is a 
$\mathcal{N}(G)$-submodule of $\ell^2(G)^n$, in particular it is a 
finitely generated $\mathcal{N}(G)$-module. The von Neumann dimension of a 
$\mathcal{N}(G)$-module is defined for any discrete group 
(\cite[Definition 1.10]{Luck}): it assigns to a nonzero finitely generated 
$\mathcal{N}(G)$-module a positive real number and satisfies the same 
basic properties that the dimension of a vector space does 
(see Theorem 1.12 in \cite{Luck}). In general it is not 
easy to compute. However, in our case the group ring 
$\mathbb{Z}[G]$ has a fraction field $\mathbb{Q}(G)$ identified with the 
field of rational functions in $m$ variables $\mathbb{Q}(t_1,\ldots,t_m)$.  
We can then define the rank of $M$ as 
the dimension of the $\mathbb{Q}(G)$-vector space $M\otimes\mathbb{Q}(G)$. 

\begin{lem}
Let $M$ be a submodule of $\mathbb{C}[G]^n$; then the von Neumann 
dimension of its closure in $\ell^2(G)^n$ is equal to its rank.
\end{lem}

\begin{proof}
This is easy to deduce from the two lemmas below: according to Lemma 
\ref{plongement} we can pick a free module $L\subset M$ such that 
$M/L$ is torsion; then the rank of $M$ is equal to that of $L$, and by 
Lemma \ref{torsdense} the completion of $L$ is equal to that of $M$. 
Thus it suffices to show the lemma for a free submodule $L$. Let 
$r$ be the rank of $L$; then $L$ is the image of an injective morphism 
from $\mathbb{C}[G]^r$ to $\mathbb{C}[G]^n$, which extends to an 
isomorphism of $\mathcal{N}(G)$-modules from $\ell^2(G)^r$ to the closure of 
$L$, so that its von Neumann dimension is $r$. 
\end{proof}

\begin{lem}
Let $M$ be a finitely generated $\mathbb{Z}[G]$-module. If $L\subset M$ is a 
maximal rank free submodule (which always exists) then the quotient $M/L$ is 
torsion.
\label{plongement}
\end{lem}

\begin{proof}
If $a_1,\ldots,a_r$ is a basis for the finite-dimensional vector space 
$M\otimes\mathbb{C}(G)$ such that $a_i\in M$ for all $i$, then 
the module $L=\sum_i\mathbb{Z}[G]a_i$ is a maximal rank free submodule. 
If $f$ is the lowest common multiple of the 
denominators of the coordinates of a generating 
family for $M$ in the basis $(a_1,\ldots,a_r)$ we have $M\subset fL$; thus 
$f$ annihilates $M/L$.
\end{proof}

\begin{lem}
The quotient of a $\mathbb{C}[G]$-module $M$ by a submodule $M'$ is torsion if 
and only $M'$ is dense in $M$.
\label{torsdense}
\end{lem}

\begin{proof}
If the quotient is torsion then for any $x\in M$ the intersection of 
$\mathbb{C}[G]x$ with $M'$ is equal to $Ix$ for some nonzero ideal $I$; thus 
it suffices to prove that an ideal is dense in $\mathbb{C}[G]$. This is true: 
any ideal $I$ contains a nonzero principal ideal $(f)$ (since we are in a 
Noetherian UFD) and the operator of $\ell^2(G)$ induced by multiplication 
by $f$ is injective, so it has dense image by \cite[Lemma 1.14]{Luck}.

If $M'$ is dense then its $\ell^2$ completion equals that of $M$ and so its 
rank, which equals its Von Neumann dimension, is the same as that 
of $M$. Now it is clear that the quotient of a module by a submodule of 
maximal rank is torsion (see the proof of the above lemma). 
\end{proof}

In the same way, given a complex $C_*,d_*$ of free finitely generated 
$\mathbb{Z}[G]$-modules we can define its reduced $\ell^2$-homology: 
the $C_i\otimes\mathbb{C}$ can be endowed with a $G$-invariant prehilbertian 
inner product (where we pick a $\mathbb{Z}[G]$-basis and its images under 
$G$ to be an orthonormal basis) and the diferentials $d_i$ extend to bounded 
operators $d_i^{(2)}$ from $C_i^{(2)}$ to 
$C_{i-1}^{(2)}$. The reduced $i$th $\ell^2$-homology module $H_i^{(2)}(C)$ 
of $C_*$ is then defined to be the $\mathcal{N}(G)$-module 
$\ker(d_i^{(2)})/\overline{\im(d_{i+1}^{(2)})}$ and the $i$th $\ell^2$-Betti 
number $b_i^{(2)}(C)$ to be its dimension. It follows from the Lemma 
above that we have $b_i^{(2)}(C)=\rg(H_i(C))$ 
(see also \cite[Lemma 1.34]{Luck}). The complex $C_*$ is thus 
$\ell^2$-acyclic if and only all $H_i(C)$ are $\mathbb{Z}[G]$-torsion modules.


\subsection{Alexander polynomials and Fuglede-Kadison determinant}

We now present the analogue for $\mathbb{Z}[G]$-modules of the order of the 
torsion subgroup for $\mathbb{Z}$-modules, and the Fuglede-Kadison 
determinant which is linked to the Mahler measure.

\subsubsection*{Alexander polynomials of $\mathbb{Z}[G]$-modules}

\label{Alexander1}

Let $M$ be a finitely generated $\mathbb{Z}[G]$-module and 
$A\in M_{m,n}(\mathbb{Z}[G])$ a presentation matrix for $M$, i.e. 
$M=\mathbb{Z}[G]^m/A\mathbb{Z}[G]^n$. We define the $l$-th elementary ideal 
of $M$ as  the ideal of $\mathbb{Z}[G]$ generated by the 
$(m-l)$-minors of $A$. This definition makes sense since a given presentation 
matrix for $M$ can be made into any other one by a sequence of elementary 
transformations that do not affect the elementary ideals (cf. 
\cite[Lemma 4.4]{Turaev}). Since $\mathbb{Z}[G]$ is a Noetherian UFD, we can 
take the greatest common divisor of an ideal; we thus put 
$\Delta_l(M)=gcd(I_l)$, and call it the $l$-th Alexander polynomial 
of the module $M$. 
The Laurent polynomial $\Delta_l(M)$ is of course defined only up to 
multiplication by a unit of $\mathbb{Z}[G]$ and it is zero for all 
$l\le \rg(M)$. It depends 
only on the $\mathbb{Z}[G]$-torsion submodule of $M$; in fact we have 
(\cite[Chapter 4]{Turaev}):
$$\Delta_{\rg(M)+l}(M)=\Delta_l(M_{tors}).$$
In particular, it follows from this equality that $\Delta_r(M)=1$ for 
a torsion-free module $M$ of rank $r$. 
We remark that $\Delta_0$ is a multiplicative invariant, i.e. if 
we have a short exact sequence 
$0\rightarrow M'\rightarrow M\rightarrow M''\rightarrow0$ of modules, 
then $\Delta_0(M)=\Delta_0(M')\Delta_0(M'')$. To prove this we can suppose 
the modules are torsion (if not both sides of the equality are zero). 
Thus $M',M''$ have presentation matrices $A',A''$ which are of maximal 
ranks $m',m''$; then $M$ has a presentation matrix of the form:
$$A=\begin{pmatrix}A'&*\\0&A''\end{pmatrix},$$
whose nonzero $(m'+m'')$-minors are products of $m'$-minors of $A'$ with 
$m''$-minors of $A''$.

Note that the Mahler measure of the Alexander polynomial $\Delta_l(M)$ is a 
well-defined positive real number since the units of $\mathbb{Z}[G]$ are the 
monomials $\pm t^v$ for $v\in\mathbb{Z}^m$, whose Mahler measure is equal to 
one.

\subsubsection*{Fuglede-Kadison determinant}

The Fuglede-Kadison determinant is defined for any $\mathcal{N}(G)$-module 
morphism, cf. \cite[Definition 3.11]{Luck}; it is an infinite-dimensional 
analogue of the $\pdet$ introduced in the next subsection. Here we only need 
it for those 
morphisms that are induced from morphisms between free $\mathbb{Z}[G]$-modules.
 It is always nonzero when the modules are finitely generated. If $A$ is 
a $(n\times n')$-matrix with coefficients in $\mathbb{Z}[G]$ we abusively 
denote by $A$ the continuous operator from $\ell^2(G)^{n'}$ to $\ell^2(G)^n$ 
that it induces, and by $\det{}_{\mathcal{N}(G)}(A)$ its Fuglede-Kadison 
determinant. We list some of its properties below for the reader's 
convenience:

\begin{itemize}
\item[(i)] If $A$ is square with nonzero determinant, then we have:
 \begin{equation}
 \det{}_{\mathcal{N}(G)}(A)=\mathcal{M}(\det A).
 \label{FKdet}
 \end{equation}
\item[(ii)](\cite[Theorem 3.14(1)]{Luck}) If 
$U\overset{f}{\rightarrow}V\overset{g}{\rightarrow}W$, $f$ has dense image
and $g$ is injective then
 $$\det{}_{\mathcal{N}(G)}(g\circ f)=\det{}_{\mathcal{N}(G)}(f)\det{}_{\mathcal{N}(G)}(g).$$
\item[(iii)](\cite[Theorem 3.14(2)]{Luck}) If $f_1$ has dense image and $f_2$ 
is injective then 
 $$\det{}_{\mathcal{N}(G)}\begin{pmatrix} f_1 & f_3\\ 0 & f_2\end{pmatrix}=\det{}_{\mathcal{N}(G)}(f_1)\det{}_{\mathcal{N}(G)}(f_2).$$
\item[(iv)](\cite[Lemma 3.15(3)]{Luck}) We have:
 $$\det{}_{\mathcal{N}(G)}(f|_{\ker(f)^{\bot}}^{\overline{\im(f)}})=\det{}_{\mathcal{N}(G)}(f)$$
\item[(v)](\cite[Lemma 3.15(4)]{Luck}) Let $f^*$ be the adjoint of $f$, then:
 $$\det{}_{\mathcal{N}(G)}(f^*)=\det{}_{\mathcal{N}(G)}(f)=\sqrt{\det{}_{\mathcal{N}(G)}(ff^*)}.$$
\end{itemize}

We prove (i), which is left as an exercise in \cite{Luck}.
Let $A=LPU$ be a Bruhat decomposition of $A$ over the field  
$\mathbb{C}(G)$, i.e. $P$ is a permutation matrix, $L$ is a 
lower triangular matrix with diagonal coefficients all equal to 1 and 
$U$ is an upper triangular matrix with diagonal coefficients 
$u_1,\ldots,u_n\in\mathbb{C}(G)$. 
Let $q$ be the lcm of the denominators of the coefficients of 
$L$ and $U$. We get 
$\det_{\mathcal{N}(G)}(q^2A)=\det_{\mathcal{N}(G)}(qLPqU)$; 
The permutation matrix $P$ is diagonalisable over $\mathbb{C}$ with modulus one
 eigenvalues, and thus $\det_{\mathcal{N}(G)}(P)=1$. From properties (ii) and 
(iii) above follow first:
$$\det{}_{\mathcal{N}(G)}(q^2Id)\det{}_{\mathcal{N}(G)}(A)=\det{}_{\mathcal{N}(G)}(qL)\det{}_{\mathcal{N}(G)}(P)\det{}_{\mathcal{N}(G)}(qU),$$
and then:
\begin{eqnarray*}
\det{}_{\mathcal{N}(G)}(A)&=&\mathcal{M}(q)^{-2n}\mathcal{M}(q^n)\mathcal{M}(q^nu_1\ldots u_n)\\
         &=&\mathcal{M}(u_1\ldots u_n)=\mathcal{M}(\det(A)).
\end{eqnarray*}


\subsection{Reidemeister torsion and $\ell^2$-torsion}

\subsubsection*{Reidemeister torsion}

Let us first recall some elementary definitions: if $f$ is a map 
between two finite-dimensional hermitian spaces $V,W$ we define 
$\pdet(f)=det(f^*f|_{\ker(f)^{\bot}})$ (the product of all nonzero 
eigenvalues of the self-adjoint map $f^*f$). A lattice of $V$ is defined to 
be a $\mathbb{Z}$-submodule $L$ of rank $\dim(V)$ spanning $V$ over 
$\mathbb{C}$; its volume is defined to be the absolute value of the 
determinant of any 
map sending an orthonormal basis of $V$ to a $\mathbb{Z}$-basis of $L$. 
If $V'$ is a subspace 
of $V$ such that $V'\cap L$ spans $V'$, then this is also true for its 
orthogonal complement $V''$ and we have
\begin{equation}
\vol(L)=\vol(L\cap V')\vol(L'')=\vol(L')\vol(L\cap V'') 
\label{ortho}
\end{equation}
where $L',L''$ are 
the orthogonal projections of $L$ on $V',V''$. 
From \eqref{ortho} it is easy to deduce the ``metric rank formula''
\begin{equation}
\vol(\im(f))=\frac{\det{}'(f)\vol(L)}{\vol(\ker f)}.
\label{rangmetr}
\end{equation}

Now let $C_*,d_*$ be a finite complex of free $\mathbb{Z}[G]$-modules; we 
pick for each $i$ a $\mathbb{Z}$-basis of $C_i$ and endow 
$C_*\otimes\mathbb{C}$ with the inner product where all these bases are 
orthonormal. The free part $H_i(C)_{free}$ of the homology is a lattice in 
$H_i(C\otimes\mathbb{C})$, which itself can be seen as the orthogonal 
of $\im(d_{i+1})$ in $\ker(d_i)\otimes\mathbb{C}$. 
Then the Reidemeister torsion defined by:
\begin{equation*}
\tau(C_*)=\prod_i \pdet(d_i)^{(-1)^i}
\end{equation*}
does not depend on the choice of the bases and is in fact given by the 
following equality:
\begin{equation}
\tau(C_*)=\prod_i \left(\frac{|H_i(C)_{tors}|}{\vol(H_i(C)_{free})}^{(-1)^{i+1}}\right).
\label{Rtorsion}
\end{equation}
The proof of this equality is a simple exercise using \eqref{rangmetr}.

\subsubsection*{$\ell^2$-torsion}

Let $C_*,d_*$ be a finite complex of free $\mathbb{Z}[G]$-modules. By 
analogy with the finite-dimensional case we put:
$$\tau^{(2)}(C)=\prod_i\fkdet(d_i)^{(-1)^i}.$$
In our particular case this can be computed:

\begin{prop}
In the situation above, suppose that $C_*$ is $\ell^2$-acyclic. Then we have:
\begin{equation}
\tau^{(2)}(C)=\prod_i \mathcal{M}\bigl(\Delta_0(H_i(C))\bigr)
\label{L2tors}
\end{equation}
\label{propl2tors}
\end{prop}

\begin{proof}
This can be deduced without too much effort from the equality 
\cite[Theorem 4.7]{Turaev}. We do not do so here, and instead prove 
a result which is also valid for non $\ell^2$-acyclic complexes (Lemma 
\ref{calctors}), following the lines of the proof of \eqref{Rtorsion}. 
The equality \eqref{L2tors} follows immediately from this result.
\end{proof}

\subsubsection*{Example: link complements}
Let $X$ be a CW-structure on the complement of a regular neighbourhood $L$ of 
a link in the three-sphere, and put $G=H_1(X)$. Let $\widehat{X}$ be the 
maximal abelian cover of $X$; then $\widehat{X}$ is a $G$-CW-complex. 
The only $H_i(\widehat{X})$ that can contain $\mathbb{Z}[G]$-torsion is
$H_1$, and it is known that in the case of a knot $K$, $H_1$ is in fact a 
torsion module. In this case we get from Proposition \ref{propl2tors} that 
$$\tau^{(2)}(\widehat{X})=\mathcal{M}(\Delta(K))$$
where $\Delta(K)$ is the Alexander polynomial of $K$. For a link with 
nonzero first Alexander polynomial $\Delta(L)$ we get in the same way 
that $\tau^{(2)}(\widehat{X})=\mathcal{M}(\Delta(K))$; if $\Delta(L)=0$ 
then the $\ell^2$-torsion
depends on the $CW$-structure; the precise dependance is computed in 
Lemma \ref{calctors}.


\section{Approximation for the Fuglede-Kadison determinant}

Let $H\subset G$ be a subgroup and $M$ a 
$\mathbb{Z}[G]$-module. We denote by $M_H$ the $\mathbb{Z}[G/H]$-module 
obtained by restricting from $G$ to $G/H$:
$$M_H=\mathbb{Z}[G/H]\otimes_{\mathbb{Z}[G]} M.$$ 
There is a $G/H$-equivariant 
isomorphism from $M_H$ to the quotient module $M/(H-1)M$. 
If $f:\mathbb{Z}[G]^n\rightarrow\mathbb{Z}[G]^{n'}$ is a linear map of 
free modules, $f_H$ is the map from $\mathbb{Z}[G/H]^n\rightarrow
\mathbb{Z}[G/H]^{n'}$ that it induces. The restricted determinant $\det'(f_H)$ 
is always taken with respect to the canonical metric on $\mathbb{C}[G/H]$, 
for which the elements of $G/H$ are an orthonormal basis. 
Let $C_*,d_*$ be a complex of $\mathbb{Z}[G]$-modules; then $C_H$ is the 
$\mathbb{Z}[G/H]$-complex with $C_{*,H}=(C_*)_H$ and $d_{*,H}=(d_*)_H$.

It is a well-known result of L\"uck (cf. \cite{Luck2} or chapter 13 of 
\cite{Luck}) that for any complex $C_*$ of free $\mathbb{Z}[G]$-modules the 
normalised Betti numbers converge to the $\ell^2$ Betti number, i.e. 
\begin{equation}
\lim_{\alpha(H)\rightarrow\infty}\frac{\dim(H_i(C_H)\otimes\mathbb{C})}{[G:H]}=b_i^{(2)}(C_*).
\label{Luck}
\end{equation}
This result is true for any residually finite group and any exhaustive 
sequence of finite index normal subgroups. It is a natural question to ask 
whether a similar result holds for torsions, i.e. does the sequence 
$\pdet(f_H)^{\frac{1}{[G:H]}}$ converge to $\fkdet(f)$. 
We shall prove the following result.

\begin{prop}
Let $A\in M_{n,n'}(\mathbb{Z}[G])$ be any matrix. Then:
$$\limsup_{\alpha(H)\rightarrow\infty}\frac{\log(\det(A_H))}{[G:H]}=\log(\fkdet(A)).$$
\label{det}
\end{prop}

More precisely, the result we prove and use in the next section is:

\begin{prop}
Let $A_1,\ldots,A_n$ be a finite collection of matrices with coefficients in 
$\mathbb{Z}[G]$; there exists a sequence $G_N$ of subgroups of finite index 
in $G$ such that $\alpha(G_N)\underset{N\rightarrow\infty}{\longrightarrow}
\infty$ and for all i we have:
\begin{itemize}
\item[(i)] $$\lim_{N\rightarrow\infty}\frac{\log\pdet(A_{i,G_N})}{[G:G_N]}=\log\fkdet(A_i)$$
\item[(ii)] $$|\dim(\ker(A_{i,G_N})\otimes\mathbb{C})-[G:G_N]\rg(\ker(A_i))|=O(\log[G:G_N]).$$
\end{itemize}
\label{determinant}
\end{prop}

Proposition \ref{det} follows from this since arguing as in the proof of 
Lemma \ref{m=1} below it is easy to see that: 
$$\limsup_{\alpha(H)\rightarrow\infty}(\det{}'(A_H)^{\frac{1}{[G:H]}})\le\det{}_{\mathcal{N}(G)}(A)$$
(this is a general fact, see Theorem 3.4(2) in \cite{Luck2}).


\subsection{Growth of Betti numbers}

\label{Betti}

We shall need precise estimates on the 
speed of convergence in \eqref{Luck}. 
For a subgroup $H\subset\mathbb{Z}^m$ we set:
$$H^{\bot}=\{(\zeta_1,\ldots,\zeta_m)\in\mathbb{T}^m,\:\forall(v_1,\ldots,v_m)\in H,\: \zeta_1^{v_1}\ldots\zeta_m^{v_m}=1\},$$
The result we want is the following (we count the multiplicities when 
computing the number of zeroes of some polynomial):

\begin{prop}
We remind the reader that $G=\mathbb{Z}^m$ and that we identify the 
group ring $\mathbb{Z}[G]$ with the ring of Laurent polynomials 
$\mathbb{Z}[t_1^{\pm1},\ldots;t_m^{\pm1}]$.
Let $M$ be a finitely generated $\mathbb{Z}[G]$-module (resp. $C_*$ a finite 
complex of finitely generated free $\mathbb{Z}[G]$-modules). Then there exists 
a polynomial $P\in\mathbb{Z}[t_1,\ldots,t_m]$ such that for 
all subgroups $H\subset G$ we have:
$$|\dim(M_H\otimes\mathbb{C})-[G:H]\rg(M)|\le|\{\zeta\in H^{\bot},\: P(\zeta)=0\}|$$
(resp.
$$|b_i(C_H)-[G:H]b_i^{(2)}(C_*)|\le|\{\zeta\in H^{\bot},\: P(\zeta)=0\}|\text{)}.$$
\label{bornes}
\end{prop}

\begin{proof} 
Let $M$ be a finitely generated $\mathbb{Z}[G]$-module; let 
$i:L\hookrightarrow M$ be 
an embedding of a free module of maximal rank in $M$ and $T$ the torsion 
module $M/L$. We see that:
$$|\dim(M_H\otimes\mathbb{C})-[G:H]\rg(M)|\le\dim(T_H\otimes\mathbb{C})+\dim(\ker(i_H)\otimes\mathbb{C}).$$
We need to bound the right-hand side; this is done in the two next lemmas.

\begin{lem}
If $T$ is a finitely generated torsion $\mathbb{Z}[G]$-module there is 
a polynomial $P$ such that we have:
$$\dim(T_H\otimes\mathbb{C})\le|\{\zeta\in H^{\bot},\: P(\zeta)=0\}|.$$
\label{bornetors}
\end{lem}

\begin{proof}
Let $a_1,\ldots,a_s$ be generators for $T$, since they are 
torsion elements there exists a nonzero 
$P\in\mathbb{Z}[t_1,\ldots,t_m]$ such that $Pa_i=0$ for all 
$i$. Then there is a surjection from $(\mathbb{C}[G]/(P))^k$ to $T$, and we 
have:
$$\dim(T_H\otimes\mathbb{C})\le k\dim(\mathbb{C}[G]/(P)\otimes\mathbb{C}[G/H]).$$
The dimension on the right is equal to $|H^{\bot}\cap\{P=0\}|$: 
we have a surjection $\mathbb{C}[G/H]\rightarrow (\mathbb{C}[G]/(P))_H$, and 
in the decomposition 
$\mathbb{C}[G/H]=\sum_{\zeta\in H^{\bot}}\mathbb{C}_{\zeta}$
(where $\mathbb{C}_{\zeta}$ is a complex line on which $G/H$ acts by the 
character associated with $\zeta$) the kernel of this surjection is 
$\bigoplus_{\zeta,P(\zeta)\not=0}\mathbb{C}_{\zeta}$. Thus 
$\dim(T_H)$ is bounded by the number of zeroes of $P^k$ (recall that we count 
multiplicities).
\end{proof}

\begin{lem}
Let $i:M_1\rightarrow M_2$ be an injective map of $\mathbb{Z}[G]$-modules. 
Then there is a polynomial $Q$ such that for all $H$ we have:
$$\dim(\ker(i_H)\otimes\mathbb{C})\le |\{\zeta\in H^{\bot},\: Q(\zeta)=0\}|.$$
\label{borneker}
\end{lem}

\begin{proof} 
We first consider the case where $M\subset\mathbb{Z}[G]^n$ is a submodule and 
$i$ its embedding. We begin by proving there is a $Q_1$ such that:
\begin{equation}
\dim(M_H\otimes\mathbb{C})\le[G:H]\rg(M)+|\{\zeta\in H^{\bot},\: Q_1(\zeta)=0\}|.
\label{(i)}
\end{equation}
To prove this take a torsion module $T=M/L$ where $L\subset M$ is a free 
maximal rank submodule: the image of $L_H\otimes\mathbb{C}$ in 
$M_H\otimes\mathbb{C}$ has dimension less than 
$\dim(L_H\otimes\mathbb{C})=\rg(M)[G:H]$ and codimension equal to 
$\dim(T_H\otimes\mathbb{C})$ so that we can take for $Q_1$ the polynomial 
associated to $T$ by Lemma \ref{bornetors}.

Now let $M'\subset\mathbb{Z}[G]^n$ be a submodule such that $M\cap M'=0$ and 
$\rg(M)+\rg(M')=n$ and $i'$ its embedding. The quotient 
$T':=\mathbb{Z}[G]^n/M\oplus M'$ is a torsion 
module and so from Lemma \ref{bornetors} we get a polynomial $Q_2$ such that 
\begin{equation}
\dim(T_H'\otimes\mathbb{C})\le|\{\zeta\in H^{\bot},\: Q_2(\zeta)=0\}|.
\label{ii}
\end{equation}
On the other hand we have an exact sequence:
$$0\rightarrow M\oplus M'\rightarrow \mathbb{Z}[G]^n\rightarrow T'\rightarrow 0$$
which becomes after tensorization with $\mathbb{Z}[G/H]$:
$$M_H\oplus M_H'\rightarrow \mathbb{Z}[G/H]^n\rightarrow T_H'\rightarrow 0.$$
We get:
$$\dim(\ker(i_H)\otimes\mathbb{C})+\dim(\ker(i_H')\otimes\mathbb{C})=\dim(T_H'\otimes\mathbb{C})+\dim(M_H\oplus M_H'\otimes\mathbb{C})-n[G:H].$$
We can apply the reasoning leading to \eqref{(i)} to $M'$ and this 
yields a polynomial $Q_3$ such that $\dim(M_H'\otimes\mathbb{C})-[G:H]\rg(M')$ 
is less than the number of zeroes of $Q_3$ on $H^\bot$; 
letting $Q=Q_1Q_2Q_3$, \eqref{(i)} and \eqref{ii} above and this imply that the 
right hand side is bounded above by the number of zeroes of $Q$ in 
$H^{\bot}$.

In the general case we have an embedding of $M_1/(M_1)_{tors}$ in a free 
module $\mathbb{Z}[G]^n$: let $i'$ be the embedding of $M_2/(M_2)_{tors}$ in 
$\mathbb{Z}[G]^n$ this yields. We get Lemma \ref{borneker} by applying the 
preceding discussion to $i'$ and Lemma \ref{bornetors} to $(M_2)_{tors}$.
\end{proof}

Let $C_*$ be a finite complex of free $\mathbb{Z}[G]$-modules; we use 
standard arguments to deduce the second point of the proposition from the 
first one. Let $n_i=\rg_{\mathbb{Z}[G]}(C_i)$; then we have 
$\rg_{\mathbb{Z}}(C_{i,H})=[G:H]n_i$, and it follows that:
\begin{align*}
b_i(C_H)&= \rg_{\mathbb{Z}}(\ker(d_{i,H}))-\rg_{\mathbb{Z}}(\im(d_{i+1,H}))\\
        &= \rg_{\mathbb{Z}}(\ker(d_{i,H}))+n_{i+1}[G:H]-\rg_{\mathbb{Z}}(\ker(d_{i+1,H})).
\end{align*}
On the other hand,
$$b_i^{(2)}(C)=\rg_{\mathbb{Z}[G]}(\ker(d_i))+n_{i+1}-\rg_{\mathbb{Z}[G]}(\ker(d_{i+1}))$$
so that we need only prove that for a matrix $A\in M_{n,n'}(\mathbb{Z}[G])$ 
the difference 
$$\rg_{\mathbb{Z}}(\ker(A_H))-[G:H]\rg_{\mathbb{Z}[G]}(\ker(A))$$
is bounded by the number of zeroes in $H^{\bot}$ of some polynomial. We know 
from the case of a module and Lemma \ref{borneker} that
$$\rg_{\mathbb{Z}}\im(A_H)=[G:H]\rg_{\mathbb{Z}[G]}(\im A)+e$$ 
where $e$ is the number of zeroes in $H^{\bot}$ of some polynomial, and since
$$\rg_{\mathbb{Z}}\im(A_H)+\rg_{\mathbb{Z}}\ker(A_H)=n'[G:H]=(\rg_{\mathbb{Z}[G]}\im(A)+\rg_{\mathbb{Z}[G]}\ker(A))[G:H]$$
we get the result we want.
\end{proof}

There is a majoration of the number of zeroes of a polynomial on the 
$m$-torus: given an affine subvariety $X\subset\mathbb{C}^m$ we denote by  
$X[H]$ the finite set $H^{\bot}\cap X$; then there exists a constant $C$ 
depending only on $X$ such that:
\begin{equation}
|X[H]|\le C\frac{[G:H]}{\alpha(H)}.
\end{equation}
We can retrieve from Proposition 
\ref{bornes} and this inequality a generalisation to the non $\ell^2$-acyclic 
case of Theorem 2.1 in \cite{Clair/Whyte} (where a proof of the above bound 
can be found).

\begin{cor}
Let $M$ be a finitely generated $\mathbb{Z}[G]$-module (resp. $C_*$ a finite 
complex of finitely generated free $\mathbb{Z}[G]$-modules). Then there exists 
a constant $C>0$ such that:
$$|\dim(M_H\otimes\mathbb{Q})-[G:H]\rg(M)|\le C\frac{[G:H]}{\alpha(H)}$$
(resp.
$$|b_i(C_H)-[G:H]b_i^{(2)}(C_*)|\le C\frac{[G:H]}{\alpha(H)}\text{)}.$$
\label{bornes2}
\end{cor}


\subsection{Convergence of Riemann sums}

We shall soon give a method to construct from a finite number of polynomials 
in $m$ variables a sequence of subgroups $G_N$ of $G$ for which the Riemann 
sums over the subsets $G_N^{\bot}$ converge to the Mahler measure, and the 
polynomials have ``few'' zeroes in $G_N^{\bot}$. We will also need a lower 
bound on the nonzero values of the polynomials at points in $G_N^{\bot}$. In 
the case where $G$ is cyclic these are two well-known facts which we recall 
below.

\begin{lem}
Let $P\in\mathbb{Z}[t]$ be a nonzero polynomial;
\begin{itemize}
\item[(i)]$\lim_{N\rightarrow\infty}\frac 1 N \sum_{\substack{\zeta^N=1\\P(\zeta)\not=0}}\log|P(\zeta)|=m(P).$
\item[(ii)] There is a $C_P>0$ such that for all $N$th roots of unity $\zeta$ 
such that $P(\zeta)\not=0$ we have: $\log|P(\zeta)|\ge-C_P\log N$.
\end{itemize}
\label{m=1}
\end{lem}

\begin{proof}
From the additivity of the logarithmic Mahler measure and the decomposition of 
$P$ into irreducible factors over $\mathbb{C}$ we see that the only case we 
need consider is that of a degree one polynomial of the form $t-a$ for some 
algebraic number $a$; if the absolute value of $a$ is not 1 then the result 
is trivial. If $|a|=1$ then it is easy to see that the superior 
limit of the Riemann sums is less than 0 (just cut off $\log|z-a|$ near 
$a$ to replace it with a continuous function with arbitrarily small mean 
value over $\mathbb{T}^1$). Jensen's formula implies that the Mahler mesure 
of $z-a$ is 0, so we have to prove that:
$$\liminf_{N\rightarrow\infty}\frac 1 N \sum_{\substack{\zeta^N=1\\\zeta\not=a}}\log|\zeta-a|\ge 0.$$
From a diophantine estimate due to Baker it is deduced in Lemma 1.11 of 
\cite{Everest/Ward} that there exists a constant $C_a>0$ such that:
\begin{equation}
\prod_{\substack{\zeta^N=1\\ \zeta\not=a}}|\zeta-a|\ge \frac{1}{N^{C_a}}.
\label{dioph}
\end{equation}
we get for all $N$ 
$$\frac 1 N \sum_{\substack{\zeta^N=1\\\zeta\not=a}}\log|\zeta-a|\ge\frac{-C_a\log N}{N}$$
and since the right-hand side converges to zero we get the desired result. 

Statement (ii) is a ``local'' version of \eqref{dioph} and is easily deduced 
from it.
\end{proof}

If $m>1$ the same argument as above shows that:
$$\limsup_{\alpha(H)\rightarrow\infty}\sum_{\substack{\zeta\in H^{\bot}\\P(\zeta)\not=0}}\log|P(\zeta)|\le m(P).$$
However, we cannot give good lower bounds for the values of $P$ at roots of 
unity and thus cannot prove the convergence. As a palliative, we use in 
the next subsection an approximation result due to Boyd and Lawton to be able 
to apply the above results to multivariate polynomials.


\subsection{Construction of the sequence $G_N$}

\begin{prop}
Let $\mathcal{P}$ be a finite collection of polynomials in $\mathbb{Z}[t_1,
\ldots,t_m]$. There exists a sequence of finite index subgroups 
$G_N\subset G=\mathbb{Z}^m$ such that 
$\alpha(G_N)\underset{N\rightarrow\infty}{\longrightarrow}\infty$ 
and for any polynomial $P\in\mathcal{P}$ we have:
\begin{gather}
\tag{i}\frac{1}{[G:G_N]}\sum_{\substack{\zeta\in G_N^{\bot}\\ P(\zeta)\not=0}}\log|P(\zeta)|\underset{N\rightarrow\infty}{\longrightarrow}m(P); \\
\tag{ii} |\{\zeta\in G_N^{\bot},\: P(\zeta)=0\}|\le\log[G:G_N]; \\
\tag{iii} \forall \zeta\in G_N^{\bot},\: \log|P(\zeta)|\ge -(\log[G:G_N])^2.
\end{gather}
\label{explicit}
\end{prop}

\begin{proof}
It is clear from the nature of the arguments below that we need only consider 
the case of a single polynomial $P\in\mathbb{Z}[t_1,\ldots,t_m]$.

For any vector $v=(v_1,\ldots,v_m)\in\mathbb{Z}^m$ we denote by $v^*$ the 
subgroup given by $\{(u_1,\ldots,u_m)\in\mathbb{Z}^m,\:\sum_i u_iv_i=0\}$. 
We put $q(v)=\alpha(v^*)$ and for a polynomial 
$P\in\mathbb{Z}[t_1,\ldots,t_m]$, $P_v(X)=P(X^{v_1},\ldots,X^{v_m})$. The 
main result of \cite{Lawton} is that:
$$m(P_v)\underset{q(v)\rightarrow\infty}{\longrightarrow}m(P).$$
For a prime $p$ we denote by $p=p_1<p_2<\ldots<p_m$ the $m$ primes following 
$p$; we put $r_i=\prod_{j\not=i}p_j$ for $i=1,\ldots,m$ and 
$r=(r_1;\ldots,r_m)$. It is easy to see that $q(r)\ge p$: if 
$\sum_i s_i r_i=0$ for some nonzero $s\in\mathbb{Z}^m$ 
and if $s_j\not=0$ we get $-ts_j r_j=p_j\sum_{i\not=j}\frac{r_i}{p_j}s_j$. 
Since the prime $p_j$ does not divide $r_j$ we must have $p_j|s_j$ and thus 
$q(r)\ge |s_j|\ge p_j\ge p$. Thus 
$m(P_r)\underset{p\rightarrow\infty}{\longrightarrow}m(P)$. For an integer 
$M\ge 0$ we define a subgroup $G_{p,M}$ of $G$ by:
$$G_{p,M}=\{v\in G,\sum_i v_ir_i\equiv0\pmod{M}\};$$
this is equal to $r^*\oplus\mathbb{Z}Mv$ for any $v$ generating a supplement 
to $r^*$, and has index $M$ in $G$. It is not hard to see that for 
$M>mp_1\ldots p_m$ we 
have $\alpha(G_{p,m})\ge p$: let $v\in G_{p,M}$, $v\not\in r^*$. Then 
$|\sum_i r_iv_i|\ge M$ and so we must have  $|v_j|>M/(m\prod_{i\not=j}p_i)$ 
for some $j$. This is bigger than $p$ for $M$ as above and the claim follows. 
Now given an $\varepsilon>0$, choose $p$ big enough so that 
$|m(P_r)-m(P)|<\varepsilon$; by Lemma \ref{m=1} there exists an integer $M_0$ 
such that for all $M\ge M_0$: 
$$\left|m(P_r)-\frac 1 M \sum_{\zeta^M=1, \, P_r(\zeta)\not=0}\log|P_r(\zeta)|\right|<\varepsilon.$$
We see that $G_{p,M}^{\bot}=\{(\zeta^{r_1},\ldots,\zeta^{r_m}),\:\zeta^M=1\}$ 
and thus we get:
$$\left|m(P)-\frac{1}{[G:G_{p,M}]}\sum_{\zeta\in G_{p,M}^{\bot},\, P(\zeta)\not=0}\log|P(\zeta)|\right|<2\varepsilon$$
for any $M\ge M_0$. 

Now the number of zeroes of $P$ in $G_{p,M}^{\bot}$ is less than the degree 
of $P_r$, which is itself less than: 
$$\max_i|r_i|\times\deg(P)\le \deg(P)p_1\ldots p_m\le \deg(P)2^{\frac{m(m+1)}{2}}p^m$$
since we have $p_i\le 2p_{i-1}$ by a well-known result of Chebyshev. Thus 
there exists an integer $M_1$ such that the number of zeroes of $P$ in 
$G_{p,M}$ is less than $\log(M)$ for all $M\ge M_1$.

From Lemma \ref{m=1} we get that for all $p$ there is a constant $C_p$ such 
that the nonzero values of $P$ at the points $\zeta\in G_{p,M}$ satisfy:
$$\log|P(\zeta)|\ge -C_p\log(M).$$
It follows that for $M\ge M_2=\exp(C_p)$ we have 
$\log|P(\zeta)|\ge -\log(M)^2$.

We can construct a sequence $G_N=G_{p_N,M_N}$ where $p_N$ 
is the $N$th prime number and $M_N$ a well-chosen positive integer: we 
take $M_N>mp_N\ldots p_{N+m-1}$ so that $\alpha(G_{p_N,M_N})\ge p_N$, 
$M_N\ge M_0,M_1,M_2$ ($M_0$ as above for, say, $\varepsilon=1/p$) so that (i), 
(ii) and (iii) hold.
\end{proof}


\subsection{Proof of Proposition \ref{determinant}}

We shall need the following elementary lemma: 

\begin{lem}
\label{det'}
Let $f$ be a linear map between Hermitian spaces $V,W$ and $V'$ 
a subspace in $V$ with codimension $d$; suppose that any singular value 
$\lambda$ of $f$ satisfies $c\le\lambda\le C$. Putting $f'=f_{|V'}$, we have:
$$d\log(c)\le\log(\det{}'(f))-\log(\det{}'(f'))\le d\log(C).$$
\end{lem}

\begin{proof}
If $\lambda_1\ge\ldots\ge\lambda_{\dim V}$ are the eigenvalues of an 
invertible positive self-adjoint endomorphism $g$ of $V$, we have that:
$$\lambda_i=\max_{F\le V,\dim(F)=i}\min_{x\in F}\frac{||gx||}{||x||};$$
so if $g=f^*f|_{\ker(f)^{\bot}}$ and $\lambda_1',\ldots,
\lambda_{\dim V-d}'$ are its eigenvalues in restriction to $V'$ (the singular 
values of $f'$), we have $\lambda_i\ge\lambda_i'\ge\lambda_{i+d}$ for all 
$i\ge 0$ and it follows that:
$$\log(\det{}'(f))-\log(\det{}'(f'))\le\sum_{i=1}^d\log(\lambda_i)\le d\log(C)$$
and that
$$\log(\det{}'(f))-\log(\det{}'(f'))\ge\sum_{i=\dim V-d}^{\dim V}\log(\lambda_i)\ge d\log(c).$$
\end{proof}

\begin{proof}[Proof of Proposition \ref{determinant}]
As in the proof of proposition \ref{explicit} it is clear that we can restrict 
ourselves to the case of a single matrix. 
The second thing to be acknowledged is that it suffices to prove the result 
for self-adjoint matrices since for any matrix $A\in M_{n,m}(\mathbb{C}[G])$ 
we have $\det_{\mathcal{N}(G)}(A)=\sqrt{\det_{\mathcal{N}(G)}(A^*A)}$. 
Thus we suppose that $A\in M_n(\mathbb{Z}[G])$ is a self-adjoint matrix.

We want to show that: 
$$\lim_{N\rightarrow\infty}\frac{\log\det'(A_{G_N})}{[G:G_N]}=\det{}_{\mathcal{N}(G)}(A)$$
for a sequence $G_N$ obtained from Proposition \ref{explicit}. 
The strategy is to use Lemma \ref{det'} above to reduce to the case when $A$ 
has a nonzero determinant. The proposition is then a straightforward 
consequence of (i) of Proposition \ref{explicit} applied to $\det(A)$. Let $L$ 
be a free maximal rank 
$\mathbb{Z}[G]$-submodule in $\im(A)$, $L'$ its image by $A$ and $h,g$ be 
isomorphisms of $\mathbb{Z}[G]^r$ ($r=\rg(A)$) with $L,L'$ so that 
the following diagram commute:
\[\begin{CD}
\mathbb{Z}[G]^r @>A'>> \mathbb{Z}[G]^r\\
@VhVV                  @VgVV\\
\mathbb{Z}[G]^n @>A>>  \mathbb{Z}[G]^n
\end{CD}\] 
with $A'$ a square matrix with nonzero determinant. Since 
$L\otimes\mathbb{C}$ is dense in $\ker(A)^{\bot}$ Lemma 3.15(3) in \cite{Luck} 
yields the equality:
$$\det{}_{\mathcal{N}(G)}(A)=\det{}_{\mathcal{N}(G)}(g)\det{}_{\mathcal{N}(G)}(A')\det{}_{\mathcal{N}(G)}(h^{-1}).$$

We write the characteristic polynomial $\det(t-A)$ as 
$\sum_{i=i_0}^n P_i(t_1,\ldots,t_m)t^i$ where $i_0$ is the rank of $\ker(A)$ 
and thus all $P_i$ are nonzero. Let $\mathcal{P}$ be the collection 
of all $P_i$ together with the nonzero coefficients of the characteristic 
polynomials of $A',\: hh^*$ and $gg^*$ and all the polynomials appearing in 
Lemma \ref{borneker} for $h$ and $g$ and the polynomials from Lemma 
\ref{bornetors} for $\im(A)/\im(h)$ and $\im(A)/\im(g)$. Let $G_N$ be the 
sequence associated to $\mathcal{P}$ by Proposition \ref{explicit}. Since 
$\det'A(\zeta)=P_i(\zeta)$ for $i=\dim\ker(A(\zeta))\ge i_0$ (and similarly 
for the other matrices) we get by (iii) of Proposition \ref{explicit} that 
for any singular value $\lambda$ of $A_{G_N},A_{G_N}',h_{G_N}$ or $g_{G_N}$ we 
have: 
\begin{equation}
\log\lambda\ge-\log[G:G_N]^2.
\label{Mino}
\end{equation}
From Lemma \ref{borneker} and (ii) of Proposition \ref{explicit} we see that 
the subspace $\ker(h_{G_N})+\ker(g_{G_N})$ of $\mathbb{C}[G/G_N]^r$ has 
dimension $\le\log[G:G_N]$. In the same way Lemma 
\ref{bornetors} yields that $\im(h_{G_N})$ and $\im(g_{G_N})=A_H\im h_{G_N}$ 
have codimension $\le\log[G:G_N]$ in $\ker(A_H)^{\bot}$. We define subspaces 
$V_1\subset\mathbb{C}[G/G_N]^r$ and $V_2\subset\mathbb{C}[G/G_N]^n$ by:
\begin{gather*} 
V_1=\ker(h_{G_N})^{\bot}\cap\ker(A_{G_N}')\cap A_{G_N}^{'-1}(\ker(g_{G_N})^{\bot})\\
V_2=h_{G_N}(V_1)
\end{gather*}
so that we have a commutative diagram:
\[\begin{CD}
    V_1     @>A_{G_N}'>> A_{G_N}'(V_1)\\
@Vh_{G_N}VV               @Vg_{G_N}VV\\
    V_2     @>A_{G_N}>>  A_{G_N}(V_2)
\end{CD}\] 
such that all arrows are isomorphisms.
We also get for $i=1,2$:
$$\dim V_i-r[G:G_N]=O(\log[G:G_N]).$$
In the sequel we use the abbreviation $u_H\sim v_H$ to mean that 
$\log(u_H)=\log(v_H)+o([G:H])$. From Lemma \ref{det'} and \eqref{Mino} above 
we see that for any linear map 
$f:\mathbb{Z}[G]^n\rightarrow\mathbb{Z}[G]^{n'}$ and any subspaces 
$F_N\subset\mathbb{C}[G/G_N]^n$ such that $\dim(F_N)=O([G:G_N]^b)$ for some 
$b<1$ we have $\det'(f_{G_N})\sim\det'(f_{G_N}|_{F_N})$. Applying this to 
$f=A'$, with Lemma \ref{det'} applied to $V=W=\mathbb{C}[G/G_N]^r$ and 
$V'=V_1$ we get:
$$\det{}'(A_{G_N}')\sim\det{}'(A_{G_N}'|_{V_1}).$$
In the same way we obtain:
\begin{eqnarray*}
\det{}'(A_{G_N}'|_{V_1})&=&\det{}'((g_{G_N}|_{A_{G_N}'V_1}^{A_{G_N}V_2})^{-1}A_{G_N}h_{G_N}|_{V_1})\\
           &=&\det(g_{G_N}|_{A_{G_N}'V_1}^{A_{G_N} V_2})^{-1}\det(A_{G_N}|_{V_2})\det(h_{G_N}|_{V_1})\\
           &\sim&\det{}'(g_{G_N})^{-1}\det{}'(A_{G_N})\det{}'(h_{G_N}),
\end{eqnarray*}
where the last line line follows from Lemma \ref{det'} applied to $f=g_{G_N}$, 
$V=\mathbb{C}[G/G_N]^r$ and $W=\im(A)_{G_N}$, and $V'=A_{G_N}'V_1$; 
then to $f=A_{G_N}$, $V=W=\im(A)_{G_N}$, $V'=V_2$ and finally to $f=h_{G_N}$, 
$V=\mathbb{C}[G/G_N]^r$, $W=\im(A)_{G_N}$ and $V'=V_1$. We conclude that:
$$\frac{\log(\det{}'(A_{G_N}))}{[G:G_N]}=\frac{\log(\det{}'(h_{G_N})^{-1}\det{}'(A_{G_N}')\det{}'(g_{G_N}))}{[G:G_N]}+o(1)$$ 
and we see that to prove the proposition it suffices to consider the case 
with nonzero determinant, since we can then apply it to the matrices $A'$, 
$hh^*$ and $gg^*$.

We suppose now that $P=\det(A)\not=0$; 
this case is proved as in \cite[Lemma 13.53]{Luck}: we decompose 
$\mathbb{C}[G/G_N]^n$ as the sum $\bigoplus_{\zeta\in G_N^{\bot}}V_{\zeta}^n$ 
where $G/G_N$ acts on the line $V_{\zeta}$ through the character induced by 
$\zeta$ so that $A_{G_N}$ is represented on $V_{\zeta}^n$ by the matrix 
$A(\zeta)$. The number of $\zeta\in G_N^{\bot}$ with $\det(A)(\zeta)=0$ is 
$\le\log[G:G_N]$ and since the singular 
values of the $A_{G_N}$ are bounded above by a constant an below by 
\eqref{Mino} we get:
\begin{eqnarray*}
\frac{1}{[G:G_N]}\log\det{}'(A_{G_N}) &=& \frac{1}{[G:G_N]}\sum_{\zeta\in G_N^{\bot}}\log\det{}'A(\zeta) \\
                                 &\sim& \frac{1}{[G:G_N]}\sum_{\zeta\in G_N^{\bot},\:\det(A)(\zeta)\not=0}\log\det(A(\zeta))
\end{eqnarray*}
Applying (i) of Proposition \ref{explicit} to the right hand side we see that 
it converges to 
$\log\det_{\mathcal{N}(G)}(A)$ as $N\rightarrow\infty$.
\end{proof}

\subsubsection*{Remark} There is a general condition under which it is known 
that the approximation of the Fuglede-Kadison determinant is true (see 
\cite[Theorem 3.4]{Luck2}). We could have checked this condition directly 
for the matrices $A_{G_N}$ using (iii) of Proposition \ref{explicit}, but 
we have chosen to give a direct proof instead.


\subsection{Growth of volumes for $\ell^2$-acyclic complexes}

We treat here the growth of the parasite term $\prod_i\vol(H_i)^{(-1)^i}$ in 
\eqref{Rtorsion}. The following proposition is an adaptation of the main 
result of \cite[section 7]{Bergeron/Venkatesh}.

\begin{prop}
Let $C_*$ be a finite complex of finitely generated $\mathbb{Z}[G]$-modules. 
Suppose that for some $i$ we have $H_i^{(2)}(C)=0$. Then for a sequence 
$G_N$ obtained from Proposition \ref{explicit} we have:
$$|\log \vol(H_i(C_{G_N})_{free})|=O(\log[G:G_N]^2).$$
\label{volumes1}
\end{prop}

\begin{proof}
For convenience we pick up the notation in \cite{Bergeron/Venkatesh} and 
denote $R_i(A)=\vol(H_i(A)_{free})$ for a finite complex of free finitely 
generated $\mathbb{Z}$-modules $A_*$. The two next lemmas are proved there:

\begin{lem}
Let a finite group $K$ act by isometries on a finite complex $A_*$ of free 
$\mathbb{Z}$-modules endowed with metrics. Let $M$ be bigger than all 
singular values of the differentials of $A_*$ an suppose that the $A_i$ are 
generated by vectors with length less than $\nu$. Suppose that all 
irreducible characters of $K$ appearing in the $K$-vector space 
$H_i(A)\otimes\mathbb{C}$ are contained in a subset $X$ and denote by $D$ 
the maximal dimension of the $X$-isotypical component of the $A_j$.
Then the following inequality holds:
$$R_i(A)\ge(M\nu|K|^5)^{-D}.$$
\label{majo2}
\end{lem}

\begin{lem}
Let $A$ be as above and $B_*$ be the dual complex 
$B_{n-j}=\hom(A_j,\mathbb{Z})$ with the dual metric. We have:
$$R_j(A_*)R_{n-j}(B_*)=1.$$
\label{dual}
\end{lem}

Now we know that the differentials of the complexes $(C_*)_H$ have their 
singular values bounded by a constant $M$ depending only on $C_*$ and that the 
$(C_j)_H$ are generated by vectors of length one. Let $X$ be the set of 
irreducible characters of $G/H$ appearing in $H_i(C_H)\otimes\mathbb{C}$ and 
$D$ the maximal dimension of the $X$-isotypical component of the $(C_j)_H$. 
Lemma \ref{majo2} yields:
$$R_i(C_H)\ge (M[G:H]^5)^{-D}.$$
Now Lemma \ref{dual} implies that:
\begin{equation}
(M[G:H]^5)^D\ge R_i(C_H)\ge (M[G:H]^5)^{-D};
\label{encvol}
\end{equation}
Since all $C_j$ are contained in $\mathbb{Z}[G]^n$ for some $n$ we see that 
for any set of characters $Y$ the dimension of the $Y$-isotypical component 
of $(C_j)_H\otimes\mathbb{C}$ is less than $n|Y|$. From Proposition 
\ref{bornes} and (ii) of Proposition \ref{explicit} we get that $D$ in 
\eqref{encvol} is an $O(\log[G:G_N])$ and thus that:
$$|\log R_i(C_{G_N})|=O(\log[G:G_N]^2).$$
\end{proof}

\begin{cor}
Let $C_*$ be a finite $\ell^2$-acyclic complex of free finitely generated 
$\mathbb{Z}[G]$-modules. There exists a sequence $G_N$ so that we have:
$$\lim_{N\rightarrow\infty}\left(\prod_i|H_i(C_{G_N})_{\mathbb{Z}-tors}|^{(-1)^i}\right)^{\frac{1}{[G:G_N]}}=\tau^{(2)}(C).$$
\end{cor}

\begin{proof}
From Proposition \ref{volumes1} we get that we can choose $G_N$ so that for 
all $i$ we have 
$$\frac{|\log \vol(H_i(C_{G_N})_{free})|}{[G:G_N]}$$
tends to 0 for all $i$. From Proposition \ref{determinant} applied to the 
$d_i$ we get that we can choose $G_N$ so that furthermore
$$\frac{\log\tau(C_{G_N})}{[G:G_N]}\xrightarrow[N\rightarrow\infty]{}\tau^{(2)}(C).$$
The corollary then follows from \eqref{Rtorsion}.
\end{proof}


\section{Proof of the main theorems}

We recall notation from the Introduction: if $X$ is a CW-complex with 
universal covering $\widetilde{X}$ and a 
surjection $\phi:\pi_1(M)\rightarrow G\cong\mathbb{Z}^m$, for any subgroup 
$H\subset G$ we denote by $X_H$ the Galois covering of $X$ given by 
$\phi^{-1}(H)\setminus\widetilde{X}$; its Galois group is $G/H$. We also 
denote the free abelian covering of $X$ given by 
$\ker(\phi)\setminus\widetilde{X}$ by $\widehat{X}$.
Theorem \ref{Main1} and Theorem \ref{Main2} are immediate consequences of the 
following result:

\begin{theo}
For all $i$, let $\Delta_j(H_i(\widehat{X}))$ be the first nonzero Alexander 
polynomial. There exists a sequence of subgroups $G_N\subset G$ with 
$\alpha(G_N)\xrightarrow[N\rightarrow\infty]{}\infty$ such that:
$$\lim_{N\rightarrow\infty}\frac{\log |H_i(X_{G_N})_{tors}|}{[G:G_N]}= m(\Delta_j(H_i(\widehat{X}))).$$
When $m=1$ we have:
$$\lim_{N\rightarrow\infty}\frac{\log |H_i(X_N)_{tors}|}{N}= m(\Delta_j(H_i(\widehat{X}))).$$
\label{approx}
\end{theo}

We begin by proving this for $m=1$, and then use Proposition \ref{explicit} 
and the Bombieri-Zannier theorem to deduce the result for $m>1$.


\subsection{When $G$ is cyclic}

Let $M$ be a finitely generated $\mathbb{Z}[t^{\pm1}]$-module of rank $r$; 
put $M_N=M/(t^N-1)M$. We use the same notation for a finite complex of free 
modules $C_*,d_*$. We want to show that:
\begin{equation}
\lim_{N\rightarrow\infty}\frac{\log|H_i(C_N)_{\mathbb{Z}-tors}|}{N}=m\bigl(\Delta_{b_i^{(2)}(C)}(H_i(C))\bigr)
\label{cyclhom}
\end{equation}
and also the corresponding limit for the module $M$:
\begin{equation}
\lim_{N\rightarrow\infty}\frac{\log|(M_N)_{\mathbb{Z}-tors}|}{N}=m\bigl(\Delta_r(M)\bigr)
\label{cyclmod}
\end{equation}
We shall deduce \eqref{cyclhom} from \eqref{cyclmod}. The latter is proved in 
three steps: first we consider torsion-free modules, then torsion modules and 
finally we combine those two to prove the general case.


\subsubsection*{Torsion-free modules}

Suppose that $M$ is torsion-free; we want to show that
\begin{equation}
\frac 1 N \log|(M_N)_{tors}|\underset{N\rightarrow\infty}{\longrightarrow} 0.
\label{tfcycl}
\end{equation}
We will in fact show that the torsion in $M_N$ is bounded independantly of 
$N$. We can pick an embedding of $M$ into a 
free module $L\cong\mathbb{Z}[G]^n$. Put $T=L/M$; then the 
$\mathbb{Z}$-torsion part of $T$ has a finite exponent: the 
$\mathbb{Z}$-torsion part of $T$ is isomorphic to $M'/M$ where 
$M'=L\cap M\otimes\mathbb{Q}$. The submodule $M'$ is finitely generated by 
$x_1,\ldots,x_k\in L$, and for each index $i$ there is an integer $n_i$ so 
that $n_ix_i\in M$. It follows that the exponent of $M'/M$,and thus that of 
$T$, divides $n_1\ldots n_k$. We can use this to bound the exponent of 
$M_N$:

\begin{lem}
For all $N$ the exponent of $(M_N)_{tors}$ divides that of the 
$\mathbb{Z}$-torsion part of $T=L/M$.
\label{bound}
\end{lem}

\begin{proof}
Suppose that $x\in M$ descends to an 
element of order exactly $n$ in $M_N$; this means that 
there exists $y\in M$, $y\not\in dM$ for any divisor $d>1$ of $n$ such that 
$nx=(t^N-1)y$. In the free module $L$ we see (by applying the Gauss lemma 
coordinate by coordinate) that $n$ must divide $y$. It follows that $n^{-1}y$ 
has order exactly $n$ in $L/M$.
\end{proof}

For the torsion-free module $M$, it follows from the above Lemma that the 
$\mathbb{Z}$-torsion in $M_N$ is always $K$-torsion with a $K$ not depending 
on $N$. Now \eqref{tfcycl} follows from the following lemma.

\begin{lem}
Let $M$ be a finitely generated $\mathbb{Z}[t^{\pm1}]$-module (resp. 
$C_*$ a finite $\mathbb{Z}[t^{\pm1}]$-complex) and $K$ a 
nonzero integer. Suppose that the $K$-torsion part of $M$ (resp. $H_i(C_*)$ is 
finite; then the $K$-torsion in $M_N$ (resp. $H_i(C_N)$) is bounded 
independantly of $N$. 
\label{ptors}
\end{lem}

\begin{proof}
Let $p_1^{k_1}\ldots p_l^{k_l}$ be the prime factorization of $K$. The 
$K$-torsion part of $M_N$ can be written as:
$$(M_N)_{K-tors}\cong\bigoplus_{i=1}^l\left(\bigoplus_{j=1}^{d_{i,N}}\mathbb{Z}/p^{\alpha_{i,j}}\mathbb{Z}\right)$$
where $d_{i,N}=\dim(M_N\otimes\mathbb{F}_{p_i})-\rg_{\mathbb{Z}}(M_N)$ and 
$\alpha_{i,j}\le k_i$. It follows immediately that
$$|(M_N)_{K-tors}|\le \sum_{i=1}^lp_i^{k_i}(\dim(M_N\otimes\mathbb{F}_{p_i})-\dim(M_N\otimes\mathbb{Q})).$$
Thus, we need only prove that for $p=p_1,\ldots,p_l$ the difference:
$$\dim(M_N\otimes\mathbb{F}_p)-\dim(M_N\otimes\mathbb{Q})$$
is bounded. This can be seen by examining the proof of Theorem 0.2(i) in 
\cite{LLS} in our very special case; we give a short direct proof here. Note 
that the proof is very similar to that of Proposition \ref{bornes}.

Suppose that $M$ is torsion-free, 
let $r=\rg_{\mathbb{Z}[t^{\pm 1}]}(M)$, and choose an 
embedding of $M$ into a free module $L$ of rank $r$. 
Since we know by Proposition \ref{bornes} that $\rg_{\mathbb{Z}}(M_N)-rN$ is 
bounded, we need to show that $\dim(M_N\otimes\mathbb{F}_p)-rN$ is 
bounded.  In the sequel we denote 
abusively $M,L$ the $\mathbb{F}_p$-modules $M,L\otimes\mathbb{F}_p$. 
There exists a $f\in\mathbb{F}_p[t]$ such that $fL\subset M$. 
We have $f_NL_N\subset M_N\subset L_N$; thus we only need to show that 
$\dim_{\mathbb{F}_p}\ker(f_N)$ is bounded. Suppose that 
$L=\mathbb{F}_p[t^{\pm 1}]$; then we can identify 
$L_N$ with the polynomials in $\mathbb{F}_p[t]$ of degree less than $N$. Let 
$\ell=\deg(f)$. Then if $h\in\mathbb{F}_p[t]$ is nonzero with $\deg(h)<N-\ell$ 
the product $fh$ is still a polynomial of degree $<N$ and thus $fh$ is not 
zero in $L_N$. It follows that $f_N$ is injective on a subspace of dimension 
$N-\ell$, so that its kernel has dimension less than $\ell$. If 
$L\cong\mathbb{F}_p[t^{\pm 1}]^n$ we can make the same argument 
coordinate by coordinate.

In general, according to the hypothesis on the $p$-torsion in $M$ we have 
$\rg_{\mathbb{F}_p[t^{\pm 1}]}M\otimes\mathbb{F}_p=r$. Thus, 
to deduce the general case from the torsion-free case we need only show that 
for a $\mathbb{F}_p[t^{\pm 1}]$-torsion module $T$ the dimension $\dim(T_N)$ 
is bounded. This is easily deduced from the result above by choosing a 
surjection from some $\bigl(\mathbb{F}_p[t^{\pm 1}]/(f)\bigr)^k$ to $T$.

The statement about homology is deduced from the statement for modules in the 
exact same way as in the proof of Proposition \ref{bornes}.
\end{proof}


\subsubsection*{Torsion modules}

Let $T$ be a finitely generated torsion module over $\mathbb{Z}[t^{\pm 1}]$; 
we want to show that:
\begin{equation}
\lim_{N\rightarrow\infty}\frac{\log|(T_N)_{tors}|}{N}=m(\Delta_0(T)).
\label{torscycl}
\end{equation}
We first return to the more general situation where $G\cong\mathbb{Z}^m$ to 
make a number of considerations. 

\begin{lem}
Let $T$ be a finitely generated $\mathbb{Z}[G]$-torsion module. Then there 
exists a finite resolution of $T$ by finitely generated free modules, 
which we write as:
\begin{equation}
0\rightarrow F_{m+2}\overset{\phi_{m+2}}{\rightarrow}\ldots\overset{\phi_{2}}{\rightarrow}F_1\overset{\phi_{1}}{\rightarrow}T\rightarrow 0
\label{reso}
\end{equation}
The complex $F_*,\phi_*$ is a $\ell^2$-acyclic complex and its 
$\ell^2$-torsion equals $\mathcal{M}(\Delta_0(T))^{-1}$. Moreover, 
for all $i>1$ there exists finitely generated torsion-free modules $M_i$ such 
that for any subgroup $H\subset G$ the homology $H_i(F_H)$ embeds into $M_{i,H}$.
\label{resolution}
\end{lem}

\begin{proof}
The ring $\mathbb{Z}[t_1,\ldots,t_m]$ is a regular ring with dimension $m+1$ 
by \cite[IV.D, Proposition 25]{Serre}; since its localizations are the same 
as those of $\mathbb{Z}[t_1^{\pm 1},\ldots,t_m^{\pm 1}]\cong\mathbb{Z}[G]$ the 
latter is also regular, with the same dimension. Now this means that all 
finitely generated modules 
over $\mathbb{Z}[G]$ have a projective resolution of length $\le m+2$. Since 
it is known that all projective $\mathbb{Z}[G]$-modules are free (cf. 
\cite[Corollary V.4.12]{Lam}) we have in fact that this holds with a free 
resolution, so that we get \eqref{reso}. The complex $F_*$ is given by the 
sequence
$$0\rightarrow F_{m+2}\rightarrow\ldots\rightarrow F_1\rightarrow 0$$
so that its homology is zero for $i>1$ and $T$ for $i=1$.
Thus it is a finite $\ell^2$-acyclic complex of free finitely generated 
$\mathbb{Z}[G]$-modules, and Proposition \ref{L2tors} allows to compute that 
its $\ell^2$ torsion equals $-m(\Delta_0(T))$.

Since 
$\im(\phi_{i,H})$ is equal to the image of $\im(\phi_i)_H=\ker(\phi_i)_H$, 
the homology $H_i(F_H)$ is equal to the quotient 
$\ker(\phi_{i,H})/\ker(\phi_i)_H$. The map $\phi_i$ embeds 
$M_i=F_i/\ker(\phi_i)$ into $F_{i-1}$; $\ker(\phi_{i,H})/\ker(\phi_i)_H$ is 
isomorphic to the kernel of the induced map $M_{i,H}\rightarrow F_{i-1,H}$, 
which is of course embedded in $M_{i,H}$.
\end{proof}

Now we return to $m=1$: it follows from the above Lemma and 
\eqref{tfcycl} that for $i>1$ we have
$$\frac 1 N \log|H_i(C_N)_{tors}|\le\frac 1 N \log|(M_{i,N})_{tors}| \underset{N\rightarrow\infty}{\longrightarrow} 0.$$
Thus \eqref{torscycl} follows from the Corollary to Proposition 
\ref{volumes1} applied to the complex $F_*,\phi_*$ since 
$H_1(F_N)\cong T_N$.


\subsubsection*{Conclusion}

Now we can prove \eqref{cyclmod}: we have the exact sequence
$$0\rightarrow M_{tors}\overset{i}{\rightarrow} M\rightarrow M'\rightarrow 0$$
where $M'=M/M_{tors}$ is torsion free. Tensoring with 
$\mathbb{Z}[\mathbb{Z}/N\mathbb{Z}]$ yields
$$(M_{tors})_N\overset{i_N}{\rightarrow} M_N\rightarrow M_N'\rightarrow 0.$$

\begin{lem}
The map $i_N:(M_{tors})_N\rightarrow M_N$ is injective for all $N$.
\label{injective}
\end{lem}

\begin{proof}
Suppose that $x\in M_{tors}$ and the image of $x$ in 
$(M_{tors})_N$ is in the kernel of $i_N$; this means that there exists a 
$y\in M$ such that $x=(t^N-1)y$. But it follows from this last equality that 
$y\in M_{tors}$, so that $x$ is already zero in $(M_{tors})_N$.
\end{proof}

Thus, for the module $M$ we have the inequality
$$1\le\frac{|(M_N)_{\mathbb{Z}-tors}|}{|(M_{tors})_N)_{\mathbb{Z}-tors}|}\le (M_N')_{\mathbb{Z}-tors}.$$
Since we have already proved the term on the right is bounded it follows that
\begin{align*}
\lim_{N\rightarrow\infty}\frac 1 N \log|(M_N)_{\mathbb{Z}-tors}| &=\lim_{N\rightarrow\infty}\frac 1 N \log|\bigl((M_{tors})_N\bigr)_{\mathbb{Z}-tors}|\\
       &=m(\Delta_0(M_{tors}))\\
       &=m(\Delta_r(M))
\end{align*}
which finishes the proof of \eqref{cyclmod}.\\

If $C_*,d_*$ is a finite complex of free finitely generated 
$\mathbb{Z}[t^{\pm 1}]$-modules the map from $\ker(d_i)_N$ to 
$\ker(d_{i,N})$ induces a map 
$j_N:H_i(C)_N\rightarrow H_i(C_N)$. 
The same proof as that of Lemma \ref{injective}\footnote{Suppose that $x\in\ker(d_i)$ is zero in $H_i(C_N)$, i.e. 
$x=y+(t^N-1)z,\: y\in\im(d_{i+1}),\: z\in(t^N-1)C_i$. We get that 
$(t^N-1)d_i(z)=0$, so that in fact $z\in\ker(d_i)$ and $x$ is already zero 
in $H_i(C)_N$.
} 
yields that 
$j_N$ is injective so that we get:
$$1\le\frac{|H_i(C_N)_{tors}|}{|(H_i(C)_N)_{tors}|}\le |\coker(j_N)_{tors}|.$$
To deal with the right-hand side we use the same trick as to prove 
\eqref{torscycl}: the cokernel of $j_N$ embeds into 
$\ker(d_{i,N})/\ker(d_i)_N$, which embeds into $M_N$ where $M=C_i/\ker(d_i)$ 
is torsion-free, so that $|\coker(j_N)_{tors}|$ is bounded. 
Summing up, we get that
$$\lim_{N\rightarrow\infty}\frac 1 N |H_i(C_N)_{tors}|=\lim_{N\rightarrow\infty}\frac 1 N |(H_i(C)_N)_{tors}|=m(\Delta_{b_i^{(2)}(C)}(H_i(C)))$$
where the last equality follows from \eqref{cyclmod} applied to $H_i(C)$.


\subsection{Torsion-free modules, $m>1$}

Our aim here is to prove the equivalent of \eqref{tfcycl} for some sequence 
$G_N$ obtained from Proposition \ref{explicit}. Recall from the proof of 
Proposition \ref{explicit} that we have defined, when $m>1$:
\begin{equation} 
G_N=r_N^*\oplus\mathbb{Z}k_Nv_N
\label{def}
\end{equation}
where $r_N,v_N\in G$ so that the scalar product $(r_N,v_N)$ is 
equal to 1 and $\alpha(r_N^*)$ tends to infinity as $N$ does. By taking 
$k_N$ large enough we can ensure that this sequence satifies the 
conclusions (i) and (ii) of Proposition \ref{determinant}, for any given 
finite collection of matrices. We want to show that given a torsion-free 
$\mathbb{Z}[G]$-module $M$ we can also choose $k_N$ so that:
\begin{equation}
\frac{1}{[G:G_N]}\log|(M_{G_N})_{tors}|\underset{N\rightarrow\infty}{\longrightarrow} 0
\label{tf}
\end{equation}
holds. We remark that $M_{r_N^*}$ is a $\mathbb{Z}[t^{\pm v_N}]$ module 
(for $v\in\mathbb{Z}^m$ we use the notation $t^v$ to denote the variable 
$t_1^{v_1}\ldots t_m^{v_m}$). 
The strategy is then to apply the results of the preceding section to the 
$\mathbb{Z}[t^{\pm v_N}]$-modules $M\otimes\mathbb{Z}[G/r_N^*]$ together with 
the following result, the proof of which is postponed until the Appendix 
\ref{loose}.

\begin{prop}
Let $M$ be a $\mathbb{Q}[G]$-module embedded in a free module 
$L$; then for $N$ big enough there exists a product of cyclotomic polynomials 
$\Phi_N$ such that the kernel of $M_{r_N^*}\rightarrow L_{r_n^*}$ is 
annihilated by $\Phi_N(t^{\pm v_N})$.
\label{aux}
\end{prop}

Let $M$ be a torsion-free $\mathbb{Z}[G]$-module which we embed into a free 
module $L$. It follows from the 
above Proposition that for $N$ big enough the kernel $T$ of the map 
$M_{r_N^*}\rightarrow L_{r_N^*}$ is annihilated by $K\Phi_N$ where $K$ is 
the exponent of the $\mathbb{Z}$-torsion in $L/M$ 
and $\Phi_N$ is a product of cyclotomic polynomials in the variable 
$t^{v_N}$. 
The submodule $T$ is the $\mathbb{Z}[t^{\pm v_N}]$-torsion part of 
$M_{r_N^*}$, so that the maps from $T_N=T/(t^{k_Nv_N}-1)T$ to $M_{G_N}$ is 
injective according to Lemma \ref{injective}. 
Let $M'$ be the 
image of $M_{r_N^*}$ in $L_{r_N^*}$, and put $M_N'=M'/(t^{k_Nv_N}-1)M'$; it 
follows from the exact sequence:
$$0\rightarrow T_N\rightarrow M_{G_N}\rightarrow M_N'\rightarrow 0$$
that the $\mathbb{Z}$-torsion in $M_{G_N}$ is less than the product 
$|(T_N)_{tors}|\times|(M_N')_{tors}|$. 
The $\mathbb{Z}[t^{\pm v_N}]$-module $T'=KT$ is annihilated by the 
polynomial $\Phi_N$, so that the order of $(T_k')_{tors}$ is 
an $o(k)$ according to \eqref{torscycl}. The torsion in $T_N$ is less than 
the $K$-torsion times the torsion in $T_{k_N}'$, and the former is bounded 
(in $k$) according to Lemma \ref{ptors}; it follows that we can 
choose $k_N$ so that $|(T_N)_{tors}|<k_N/N$ (say). Since 
$M_N'$ is $\mathbb{Z}[t^{\pm v_N}]$-torsion-free \eqref{tfcycl} allows us 
to choose $k_N$ so that $|(M_N')_{tors}|<k_N/N$. The conclusion 
\eqref{tf} follows (recall that $k_N=[G:G_N]$).


\subsection{Conclusion}

We are now in position to prove the following theorem, from which Theorem 
\ref{approx} follows immediately.

\begin{theo}
\begin{itemize}
  \item[(i)] Let $M$ be a finitely generated module of rank $r$; then:
    \begin{equation}
     \lim_{N\rightarrow\infty}\frac{\log|(M_{G_N})_{tors}|}{[G:G_N]} = m(\Delta_r(M)).
     \label{modules}
    \end{equation}
  \item[(ii)] Let $C_*,d_*$ be a complex of finitely generated free modules, 
      $r=\rg H_i(C)$; we have:   
    \begin{equation}
     \lim_{N\rightarrow\infty}\frac{\log|H_i(C_{G_N})_{tors}|}{[G:G_N]} = m(\Delta_r(H_i(C)))
     \label{complexes}
    \end{equation}
\end{itemize}
\label{torsion3}
\end{theo}

\begin{proof}
We use the same scheme of proof as in the cyclic case: we begin by proving 
\eqref{modules} for a torsion module, then deduce the case of a finitely 
generated module, and finally use it to prove \eqref{complexes}.

\subsubsection*{Modules}
First of all, the proof of \eqref{modules} does not change from the case 
$m=1$: Lemma \ref{resolution}, the corollary to Proposition \ref{volumes1} 
and \eqref{tf} together 
imply that for a torsion module $T$ we can choose $k_N$ in \eqref{def} so that
\begin{equation}
\frac{1}{[G:G_N]}\log|(T_{G_N})_{tors}|\underset{N\rightarrow\infty}{\longrightarrow}m(\Delta_0(T)).
\label{tors}
\end{equation}

Now let $M$ be any finitely generated $\mathbb{Z}[G]$-module; we have an 
exact sequence:
$$0\rightarrow M_{tors}\rightarrow M\rightarrow M'\rightarrow 0$$
where $M'=M/M_{tors}$. According to \eqref{tf} we can choose the $k_N$ so that 
$M_{G_N}'$ has a negligible torsion, so that to deduce \eqref{modules} from 
\eqref{tors} we need only show that:
\begin{equation}
\frac{\log|\ker\bigl((M_{tors})_{G_N}\rightarrow M_{G_N}\bigr)_{\mathbb{Z}-tors}|}{[G:G_N]}\underset{N\rightarrow\infty}{\longrightarrow} 0.
\label{noyau}
\end{equation}
The following result is deduced from Proposition \ref{aux} in Appendix 
\ref{loose}.

\begin{prop}
For $N$ big enough the kernel of the map 
$(M_{tors})_{r_N^*}\rightarrow M_{r_N^*}$ is annihilated by $K\Phi_N$ where 
$\Phi_N$ is a cyclotomic polynomial and $K$ an integer (not depending on $N$).
\label{comparison}
\end{prop}

We use this and an argument similar to that used to prove Lemma 
\ref{injective} 
to deduce \eqref{noyau}. Suppose that $x\in M_{tors}$ descends to 
an element in $\ker\bigl((M_{tors})_{G_N}\rightarrow M_{G_N}\bigr)$; then 
there exists $y\in M, z\in(r_N^*-1)M$ such that $x=(t^{k_Nv_N}-1)y+z$. It 
follows that $y\in M_{tors}+(r_N^*-1)M$, so that we can write 
$x=(t^{k_Nv_N}-1)y'+z'$ where $y'\in M_{tors}$ and $z'\in(r_N^*-1)M$ 
It follows that $x$ and $z'$ have the same image in $(M_{tors})_{G_N}$, 
so that the kernel $\ker\bigl((M_{tors})_{G_N}\rightarrow M_{G_N}\bigr)$ 
is the image of $T=\ker\bigl((M_{tors})_{r_N^*}\rightarrow M_{r_N^*}\bigr)$. 
According to Proposition \ref{comparison}, \eqref{ptors} and 
\eqref{torscycl} we can choose $k_N$ so that $(T_{k_N})_{tors}$ has 
order less than $k_N/N$. This implies that the torsion in the kernel 
of the map $(M_{tors})_{G_N}\rightarrow M_{G_N}$ also has order less 
than $k_N/N$ and \eqref{noyau} follows.

\subsubsection*{Complexes}
Let $N$ be fixed. We compare the growth of the $\mathbb{Z}$-torsion of 
$H_i(C)_{r_N^*}/(t^{kv_N}-1)H_i(C)_{r_N^*}$ and 
$H_i(C_{r_N^*}/(t^{kv_N}-1)C_{r_N^*})$; more precisely we show that:
\begin{equation}
\lim_{k\rightarrow\infty}\frac 1 k \bigl(\log|(H_i(C)_{r_N^*}/(t^{kv_N}-1)H_i(C)_{r_N^*})_{\mathbb{Z}-tors}|-\log|(H_i(C_{r_N^*}/(t^{kv_N}-1)C_{r_N^*}))_{\mathbb{Z}-tors}|\bigr)=0.
\label{gy}
\end{equation}
Let $\phi$ denote the map from $(H_i(C)_{r_N^*})_{tors}$ to 
$H_i(C_{r_N^*})_{tors}$; we prove that there exists a cyclotomic polynomial 
$\Phi_N$ which annihilates $\ker(\phi)$. Let 
$x\in\ker(d_i)$ map to an element in $\ker(\phi)$; we can write 
$x=y+z$ where $z\in\ker(d_i)\cap(r_N^*-1)C_i$. By proposition 
\ref{aux}, since there is no $\mathbb{Z}$-torsion in $C_i/\ker(d_i)$, we get 
that there exists $\Phi_N$ such that $\Phi_N z\in(r_N^*-1)\ker(d_i)$, so that 
$\Phi_N x$ maps to 0 in $(H_i(C)_{r_N^*})_{tors}$. As for 
$\coker(\phi)$, applying Proposition \ref{aux} to the embedding of 
$C_i/\ker(d_i)$ into $C_{i-1}$ we get that there is an integer $K$ independant 
of $N$ and a cyclotomic polynomial $\Phi_N'$ such that 
$K\Phi_N'\coker(\phi)=0$. It follows from these remarks and \eqref{cyclmod}, 
\eqref{cyclhom} that:
\begin{multline*}
0\le \lim_{k\rightarrow\infty}\frac 1 k \bigl(\log(H_i(C_{r_N^*}/(t^{kv_N}-1)C_{r_N^*}))-\log(H_i(C)_{r_N^*}/(t^{kv_N}-1)H_i(C)_{r_N^*})\bigr)\\
   \le \lim_{k\rightarrow\infty}\frac 1 k \log|(H_i(C_{r_N^*}/(t^{kv_N}-1)C_{r_N^*}))_{K-tors}|.
\end{multline*}
It follows from Lemma \ref{ptors} that the limit on the right is zero, proving 
\eqref{gy}. By choosing $k_N$ big enough we can thus suppose that 
$$\lim_{N\rightarrow\infty}\frac{\log|(H_i(C)_{G_N})_{\mathbb{Z}-tors}|-\log|H_i(C_{G_N})_{\mathbb{Z}-tors}|}{[G:G_N]}=0,$$
so that we deduce \eqref{complexes} from \eqref{modules} applied to 
$H_i(C)$.

\end{proof}


\section{Cyclic coverings of 3--manifolds}

In the case $m=1$, Theorem \ref{approx} yields that for any finite 
CW-complex $X$ with an infinite cyclic covering $\widehat{X}$ and 
$X_N=N\mathbb{Z}\setminus\widehat{X}$ we have:
$$\frac{\log |H_i(X_N)_{tors}|}{N}\underset{N\rightarrow\infty}{\longrightarrow}\Delta_{b_i^{(2)}(\widehat{X})}(H_i(\widehat{X})).$$
Letting $X$ be a cell structure on a compact n--manifold we get 
Theorem \ref{Main2}.

\subsection{A presentation for the Alexander module}

It would be nice to have an explicit formula for the Alexander polynomial 
$\Delta_i(\widehat{M})$. The simplest case is that of a manifold fibering over 
the circle, which is $\ell^2$-acyclic and thus already treated in 
\cite[Corollary 7.7]{Bergeron/Venkatesh}: if $\varphi$ is the diffeomorphism 
gluing the fiber $S$ and $\varphi_*$ its action on $H_1(S)$ we have 
$\Delta_0(\widehat{M})=\det(1-t\varphi_*)$. In the general case, the  
covering $\widehat{M}$ is obtained from an epimorphism 
$\pi_1(M)\rightarrow\mathbb{Z}$, and all such morphisms are given by the 
intersection form with an embedded incompressible two-sided nonseparating 
surface in $M$: see the proof of Lemma 6.6 in \cite{Hempel}. 
From here on we suppose we have a closed 3-manifold $M$ containing an embedded 
incompressible two-sided nonseparating surface $S$. Let $\widehat{M}$ be the 
infinite cyclic covering of $M$ associated with $S$. Let 
$U\cong S\times(-1,1)$ be a tubular neighbourhood of $S$ in $M$ and let 
$M'=M-U$. We denote by $S_{\pm}$ the connected components of $\partial M'$ 
corresponding to $S\times\{\pm 1\}$ in $\overline{U}$. The fundamental 
group of $M$ is then given by an HNN-extension: we consider 
$\pi_1(S_{\pm})$ as subgroups of $\pi_1(M')$, and there is an isomorphism 
$\alpha:\pi_1(S_+)\rightarrow\pi_1(S_-)$ such that
$$\pi_1(M)=\langle \pi_1(M'),t|\: \forall g\in\pi_1(S_+),\: tgt^{-1}=\alpha(g)\rangle.$$
Let us denote by $\alpha_*$ the induced map $H_1(S_+)\rightarrow H_1(S_-)$. 
We were not able 
to get a formula as explicit as that above, but only to obtain a presentation 
of $H_1(\widehat{M})$.

We can give an explicit construction of $\widehat{M}$ using $S$: 
$\widehat{M}$ is diffeomorphic to the manifold obtained from 
$M'\times\mathbb{Z}$ by identifying $S_+\times\{n\}$ with $S_-\times\{n+1\}$ 
using $\alpha$. Let $i$ be the embedding of $S$ in $M'$ corresponding to 
$S_+$. We know that $H_0(\widehat{M})=0$ and the map from $H_1(S_+)$ to 
$H_1(M')$ is injective; thus the homology long exact sequence coming from:
$$S\times\mathbb{Z}\underset{(1-t\alpha)i}{\longrightarrow}M'\times\mathbb{Z}\rightarrow\widehat{M}$$
yields the short exact sequence:
\begin{equation}
\mathbb{Z}[t^{\pm1}]\otimes H_1(S)\underset{(1-t\alpha_*)i_*}{\longrightarrow}\mathbb{Z}[t^{\pm1}]\otimes H_1(M')\rightarrow H_1(\widehat{M})\rightarrow 0.
\label{prez}
\end{equation}
which is a presentation of $H_1(\widehat{M})$. However, in all 
generality it seems rather 
difficult to compute the minors of the matrix $(1-t\alpha_*)i_*$ due to 
the unknown nature of $i_*$.


\subsection{Examples}

\subsubsection*{Computation for rational homology cylinders}
Suppose that $M',i_{\pm}$ is a rational homology cylinder, that is, the 
homology map $i_*$ induces an isomorphism between rational homologies. 
Then we can compute the Alexander polynomial, generalising a result 
for knot complements:

\begin{lem}
In the situation described above, $\widehat{M}$ is $\ell^2$-acyclic and 
$\Delta_0(H_1(\widehat{M})=\tau(M',S_+)\det(1-t\alpha_*)$ where 
$\tau(M',S_+)$ denotes the order of $H_1(M')/i_*H_1(S)$.
\end{lem}

\begin{proof}
This is an immediate consequence of \eqref{prez} since in this case 
the map $(1-t\alpha_*)\circ i_*$ is represented by a square matrix 
with nonzero determinant.
\end{proof}

In the case where $M$ is a knot exterior we retrieve the result of 
\cite[Theorem 6.5]{Lickorish}. Note also that by the Dehn-Nielsen 
theorem $\alpha$ corresponds to an element $\phi$ in the mapping class group 
of $S$ and we have $\det(1-t\alpha_*)=\det(1-t\phi_*)$, so that we also 
retrieve the case where $M$ is fibered.\\

Given an homology cylinder $M',i_{\pm}$ (i.e. $i_{\pm}$ are 
embeddings of $S$ in $\partial M'$ which induce an isomorphism in 
homology), \cite[Section 4.3]{KFC} construct by surgery on the interior 
another homology cylinder. Their construction does not affect the 
Alexander polynomial of the infinite cyclic covering of the 3-manifold 
$M=M'/\{\forall x\in S,\:i_+(x)=i_-(x)\}$ dual to $S$.

\subsubsection*{Positive $\ell^2$-Betti number}
There indeed are 
cases where the infinite cyclic covering is not $\ell^2$-acyclic. 
A somewhat artificial example is obtained as follows:  
suppose that $M_0$ is a three-manifold having an infinite cyclic covering 
$\widehat{M}_0$ and that $N$ is another three-manifold with positive $b_1$. 
Then the connected sum $M=M_0\sharp N$ has an infinite cyclic covering 
$\widehat{M}$ difeomorphic to $\widehat{M}_0\sharp(N\times\mathbb{Z})$ 
(if $M_0'$ is $M_0$ cut along a surface dual to the covering, it is obtained 
by attaching a copy of $N$ to each copy of $M_0'$ in $\widehat{M}_0$). 
Thus we have
$$H_1(\widehat{M})\cong H_1(\widehat{M}_0)\oplus(\mathbb{Z}[t^{\pm1}]\otimes H_1(N))$$
and it follows that $b_1^{(2)}(\widehat{M})\ge b_1(N)$. Note that we also 
have $\Delta_{b_1(N)+i}(H_1(\widehat{M}))=\Delta_i(\widehat{M}_0)$.\\

A more interesting example is given by manifolds 
with large fundamental group. Recall that a group is said to be large when 
it has a finite-index subgroup surjecting onto a noncyclic free group. 
Suppose that $\pi_1(M)$ surjects onto the free product 
$\mathbb{Z}*\mathbb{Z}$; then we can take the surjection from $\pi_1(M)$ 
to $\mathbb{Z}$ obtained by projection onto the second free factor. 
Then the associated infinite cyclic covering $\widehat{M}$ has 
$b_1^{(2)}(\widehat{M})\ge 1$. Indeed, the kernel of the map 
$\mathbb{Z}*\mathbb{Z}\rightarrow\mathbb{Z}/N\mathbb{Z}$ is the 
subgroup freely generated by $a^N,b,aba^{-1},\ldots,a^{N-1}ba^{1-N}$, 
so that $\pi_1(M_N)$ surjects onto a free group of rank $N$ and thus 
we have $b_1(M_N)=b_1(\pi_1(M_N))\ge N$ (any three-manifold with infinite 
fundamental group is aspherical) so that $b_1^{(2)}(\widehat{M})\ge 1$ 
by L\"uck's theorem.

Lots of arithmetic three-manifolds are known to have a large 
fundamental group, and it is expected that in fact all hyperbolic 
three-manifolds have a large $\pi_1$
\footnote{Ian Agol recently announced that he had a proof of this conjecture, following a previous work of Dani Wise on the subject.}. 
See for example \cite{Lackenby} for 
recent progress on this and its links with other conjectures in 
3-manifold topology.


\appendix

\section{Loose ends}

\label{loose}

We shall prove here the results used in the proof of Theorem \ref{torsion3}.
Recall that, for a subgroup $H\subset\mathbb{Z}^{m'}$, we denote by $H^*$ the 
subgroup of vectors orthogonal to $H$ in $\mathbb{Z}^{m'}$; $\alpha(H)$ is the 
smallest length of a nonzero vector in $H$. We introduce notation we shall use 
throughout this Appendix: if $A$ is an $m\times m'$ 
matrix with coefficients in $\mathbb{Z}$ we denote by $H_A$ the 
subgroup of $\mathbb{Z}^{m'}$ spanned by its lines. We will 
continue to denote by $A$ the induced morphism from 
$\mathbb{Z}[\mathbb{Z}^{m'}]$ to $\mathbb{Z}[\mathbb{Z}^m]$.
Recall that we identify the group ring $\mathbb{Z}[\mathbb{Z}^{m'}]$ with the 
Laurent polynomials $\mathbb{Z}[t_1^{\pm 1},\ldots,t_{m'}^{\pm 1}]$; we 
denote this identification by $v\mapsto t^v$.


\subsection{The Bombieri-Schinzel-Zannier theorem}

The following result is an immediate generalisation of a theorem by 
Bombieri and Zannier (improving on a previous result by Schinzel):

\begin{theo}
Let $f_1,\ldots,f_k\in\mathbb{Z}[\mathbb{Z}^{m'}]$ be coprime polynomials. 
Then there exists a constant $c=c(f_1,\ldots,f_m)$ such that if the gcd of 
$Af_1,\ldots,Af_k$ is not equal to a cyclotomic polynomial then 
$\alpha(H_A^*)\le c$.
\label{BZ}
\end{theo}

\begin{proof}
This follows from an easy induction using the theorem as stated in 
\cite[Theorem 45 and page 517]{Schinzel}; see also \cite[Proposition 5.2]{Le}.
\end{proof}

In the sequel we shall use the term ``cyclotomic polynomial'' to design 
multivariate generalised cyclotomic polynomials , i.e. the evaluations of a 
cyclotomic polynomial at a monomial, as well as their finite products.


\subsection{Proof of Proposition \ref{aux}}

In the sequel we denote by $A^n$ the map induced by $A$ from 
$\mathbb{Z}[\mathbb{Z}^{m'}]^n$, i.e. 
$A^n(f_1,\ldots,f_n)=(Af_1,\ldots,Af_N)$.

\begin{prop}
Let $M$ be a submodule in $\mathbb{Q}[\mathbb{Z}^{m'}]^n$; 
there exists a constant $c_M$ such that for 
all $A$ with $\alpha(H_A^*)>c_M$, there exists a cyclotomic polynomial 
$\Phi_A$ so that: 
$$\Phi_A\ker(A^n)\cap M\subset\ker(A)M.$$
\end{prop}

\begin{proof}
In this proof we denote the group ring $\mathbb{Q}[\mathbb{Z}^{m'}]$ by 
$R$, and $\mathbb{Q}[A\mathbb{Z}^{m'}]$ by $R'$. It suffices to prove the 
result when the image of $A$ is of rank one (i.e. $A$ is a line matrix).

Suppose first that $M=fR$ is a nonzero principal ideal. For 
$\alpha(H)$ big enough the support of $f$ injects into $\mathbb{Z}^{m'}/H$, 
in particular for $\alpha(H_A^*)$ big enough $Af\not=0$. The ideal 
$\ker(A)$ is prime (because $R/\ker(A)\cong R'$ is 
a domain), so that if $fg\in\ker(A)$ then $g\in\ker(A)$. This means 
that $fR\cap\ker(A)=f\ker(A)$.

We next take $M$ to be an ideal in $R$, and prove the result by induction 
on the number of generators. First we reduce to the case where $M$ is not 
contained in any principal ideal. Suppose that $f$ is the gcd of $M$; 
then we can write $M=fM'$ where $M'$ is not contained in any principal 
ideal. We have 
$\ker(A)\cap fM'=f\ker(A)\cap fM'$ by the principal case, so that if 
the result is true for the embedding of $M'$ in $fR$ it is also true for $M$.

Thus, let $f_1,\ldots,f_k$ be a minimal generating family for $M$ with no 
common factor.  Suppose that 
$\alpha(H_A^*)>c=c(f_1,\ldots,f_k)$ with 
$c$ the constant from Theorem \ref{BZ}; then the gcd of $Af_1,\ldots,Af_k$ 
must be equal to a cyclotomic polynomial $\Phi\in AR$ (we identify it with 
its evaluation at the monomial $t^{v_A}$ 
where $v_A$ generates a supplementary for $H_A^*$). 
Let $h=\sum_i h_if_i\in\ker(A)\cap M$; we get that: 
$$-Ah_1Af_1=\sum_{i\ge 2}Ah_iAf_i$$ 
and so the gcd of $Af_2,\ldots,Af_k$ divides $\Phi Ah_1$. We can write this 
gcd 
as $\sum_{i\ge 2}Ag_iAf_i$ for some $g_i\in R$ since 
$R'$ is a principal ring (because 
$A\mathbb{Z}^{m'}$ is a cyclic group). Thus we get 
$\Phi Ah_1=\sum_{i\ge 2}Ah_i'Af_i$  
for some $h_i'\in R$. It follows that we have:
$$\Phi h_1=\sum_{i\ge 2}h_i'f_i+h_1'$$
for some $h_1'\in\ker(A)$, and from that:
$$\Phi h=\Phi h_1'f_1+\sum_{i\ge 2}(h_i+h_i'f_1)f_i.$$
By the induction hypothesis, 
if $\alpha(H_A^*)>c_{M'}$ (with $M'=(f_2,\ldots,f_k)$) then 
$\sum_{i\ge 2}(h_i+h_i'f_1)f_i\in \Phi'\ker(A)M'$, so we get the result 
with $c_M=\max(c,c_{M'})$.

Suppose now that $\rg(M)=1$ and $M$ is embedded in $R^n$; up to passing to 
a larger free module (in $R\otimes\mathbb{Q}(\mathbb{Z}^{m'})^n$) we may 
suppose that $M$ is contained in a copy of $R$ that is a direct factor 
of $R^n$; then we can apply the above arguments to this embedding 
of $M$ in $R$ and get the result.

Now let us prove the general result by induction on the rank. Let $M$ have 
rank $k>1$ and be embedded in $R^n$; let $l<n$ so that the intersection 
$M_0=M\cap(R^l\times 0)$ has rank one. Put $M_1=M/M_0$; this is torsion-free 
and we have a commutative diagram:
\[\begin{CD}
0@>>>    M_0      @>>>   M   @>>>   M_1      @>>>0\\
&  &   @VVV            @VVV         @VVV  \\
0@>>> R^l\times 0 @>>> R^n   @>>>  R^{n-l}   @>>>0         
\end{CD}\]
Suppose that $x\in\ker(A^n)\cap M$; its projection $x_1$ in $M_1$ lies in 
$\ker(A^{n-l})\cap M_1$ and, by the induction hypothesis, there exists 
a cyclotomic $\Phi$ such that $\Phi x_1\in\ker(A)M_1$. 
It follows that $\Phi x=y+x_0$ for some $y\in\ker(A)M$ and 
$x_0\in \ker(A^n)\cap M_0$. By the induction hypothesis we get that there 
exists a cyclotomic $\Phi'$ such that 
$\Phi'x_0\in\ker(A)M_0$. Finally, we have that 
$\Phi\Phi' x\in\ker(A)M$.
\end{proof}

\begin{proof}[Proof of Proposition \ref{aux}]
The kernel $\ker(A)$ is easily seen to be equal to $(H_A^*-1)R$; indeed, 
the groups $A\mathbb{Z}^{m'}$ and $\mathbb{Z}^{m'}/H_A^*$ are isomorphic via 
$A$ so that the kernels of the maps 
$\mathbb{Z}[\mathbb{Z}^{m'}]\rightarrow\mathbb{Z}[A\mathbb{Z}^{m'}]$ and 
$\mathbb{Z}[\mathbb{Z}^{m'}]\rightarrow\mathbb{Z}[\mathbb{Z}^{m'}/H_A^*]$ 
must be equal. The first is $\ker(A)$ and the second 
$(H_A^*-1)\mathbb{Z}[\mathbb{Z}^{m'}]$.
Thus if $M\subset L$ is a submodule of a free module, the above Proposition 
yields that there exists a cyclotomic $\Phi_A$ 
such that $\Phi_A$ annihilates the kernel of the map
$$M_{H_A^*}=M/(H_A^*-1)M\rightarrow L/(H_A^*-1)=L_{H_A^*}.$$
Proposition \ref{aux} follows by taking $A$ to be the line matrix 
$r_N$ since $\alpha(r_N^*)$ tends to infinity.
\end{proof}


\subsection{Proof of Proposition \ref{comparison}}

We will in fact prove Proposition \ref{comparison} for any sequence of 
subgroups of finite index which all contain a direct factor and whose smallest 
length tends to infinity. The method is to use induction on the rank. Here is 
the result that allows this:

\begin{lem}
Let $H$ be a subgroup of $\mathbb{Z}^{m'}$ which contains a direct factor 
of $\mathbb{Z}^{m'}$ of rank $m'-1$, and such that $\alpha(G)>c$. Then there 
exists subgroups $H=H_1\supset H_2\supset\ldots\supset H_{m'}=0$ such that all 
inclusions are of corank one and all $H_i,\: i\ge 2$ are direct factors. 
Moreover, they all satisfy $\alpha(H_i^*)\ge c$.
\end{lem}

\begin{proof}
Let $H_2$ be a direct factor contained in $H$. It has rank $m'-1$, so the 
intersection $H_2\cap(0\times\mathbb{Z}^2)$ is not empty. Let $v=(0,a,b)$ be 
in this intersection; we may suppose that $gcd(a,b)=1$ since $H_2$ is a direct 
factor, and then we have $\min(|a|,|b|)\ge\alpha(H)>c$. Thus we see 
that $\alpha(v^*)>c$. Now we only have to choose a basis $\{v_1=v, v_2,\ldots,
v_{m'-1}\}$ of $H$ and put $H_i=\sum_{j=1}^i\mathbb{Z}v_i$; the $H_i$ 
all contain $v$ and so $\alpha(H_i^*)\ge\alpha(v^*)>c$, and it is clear that 
they satisfy the other conclusions.
\end{proof}

In consequence, Proposition \ref{comparison} follows from:

\begin{prop}
Let $M$ be a finitely generated module over $R$. 
There exists a constant $c$ depending on $M$ such that the following holds. 
If $H$ is a subgroup of  $\mathbb{Z}^{m'}$ of corank $\ge 1$ having a sequence 
of subgroups $H_2\supset\ldots\supset H_{m'}$ such that all inclusions are 
of corank one, all $H_i$ are direct factors and they all satisfy 
$\alpha(H_i^*)>c$, then the kernel of the map $(M_{tors})_H\rightarrow M_H$ 
is annihilated by a cyclotomic polynomial.
\end{prop}

\begin{proof}
If $H$ has rank one, we proved that the kernel is 0 in Lemma \ref{injective}. 
We need the following result to carry out the induction step:

\begin{lem}
There is a constant $c$ so that for $\alpha(H_A^*)>c$ there exists a 
cyclotomic polynomial $\Phi$ such that the 
map $M_{R-tors}\rightarrow \Phi_A(AM)_{R'-tors}$ 
is surjective .
\end{lem}

\begin{proof}
We put $M'=M/M_{tors}$, embed $M'$ in a free module $L$ and denote by $\pi$ 
the map $M\rightarrow L$ this yields. Suppose that $x\in M$ maps to 
a $R'$-torsion element in $AM$; this means that there exists a $f\in R$ with 
$Af\not= 0$ such that $\pi(fx)\in\ker(A)\pi(M)$. Let $c$ be the constant for 
$M'\subset L$ from Proposition \ref{aux}; we get that for $\alpha(H_A^*)>c$ 
there exists a cyclotomic $\Phi$ such that $\Phi\pi(fx)\in\ker(A)M'$. 
It follows that $\Phi x=y+x'$ where $y\in M_{R-tors}=\ker(\pi)$ and 
$x'\in\ker(A)M$, and thus that $A\Phi Ax$ is in the image of $M_{R-tors}$ 
in $AM$.
\end{proof}

Write $H=\mathbb{Z}v\oplus H_2$ where $H_2$ is a direct factor of 
$\mathbb{Z}^{m'}$. Let $x\in M_{tors}$, suppose that 
$x\in (H-1)M$. Then we can write $x=(t^v-1)x'+z$ where $z\in(H_2-1)M$. 
Let $A$ be a matrix with $\ker(A)=H_2$ (this is possible since $H_2$ is 
a direct factor); we get that $Ax=(t^{Av}-1)Ax'$, so that $x'$ is 
$R'$-torsion in $AM'$ (we have to suppose that the annihilator of $x$ is 
not contained in $(H_2-1)R$, but this is true for $\alpha(H_2)$ big enough). 
By the Lemma above we get that $\Phi x'=y+z'$ where 
$z'\in\ker(A)M=(H-1)M,\: y\in M_{tors}$. It follows that 
$\Phi x=(t^v-1)y+z+z'$, and we get that $\Phi'(z+z')\in(H-1)M_{tors}$ for 
some cyclotomic $\Phi'$, by the induction hypothesis.
\end{proof}


\section{The $\ell^2$-volume and $\ell^2$-torsion}

\label{l2volume}

We introduce here a new $\ell^2$-invariant, the aim of which is to be able to 
mimic the arguments leading to \eqref{Rtorsion}. In particular, we deduce 
\eqref{L2tors} as a particular case.

\subsection{Definition}
We want to define a $\ell^2$-volume for submodules of $\mathbb{Z}[G]^n$ 
by analogy with the lattices in a finite-dimensional Hermitian space. 
It is quite natural to do this for free submodules: if 
$L\subset\mathbb{Z}[G]^n$ is free of rank $r$, then 
all morphisms $\mathbb{Z}[G]^r\rightarrow\mathbb{Z}[G]^n$ with image $L$ have 
the same Fuglede-Kadison determinant since if we fix one such morphism $f$, 
then any other one can be written $f\circ u$ for a change of basis $u$ in 
$\mathbb{Z}[G]^r$, and all change of basis morphisms for free 
$\mathbb{Z}[G]$-modules have a unit for determinant and thus their 
Fuglede-Kadison determinant is equal to one by equation \eqref{FKdet}. Thus, 
choosing a matrix $A\in M_{n,r}(\mathbb{Z}[G])$ with image $L$ we can put, 
as in the finite-dimensional case:
$$\vol^{(2)}(L)=\fkdet(A)=\mathcal{M}(\det(A^*A)).$$
When $M$ is a submodule of $\mathbb{Z}[G]^n$ with rank $r$, we can pick 
a free submodule of maximal rank (lemma \ref{plongement}). The quotient 
$M/L$ is torsion by lemma \ref{torsdense}, so it has nonzero first 
Alexander polynomial. In the finite dimensional case, if we have a finite 
index sumodule $V'$ of a module $V$ then for any metric on $M$ 
the equality $\vol(V)=[V:V']\vol(V')$ holds. Since we are interested in 
approximation problems the analogue of the index we shall consider for the 
maximal rank submodule $L\subset M$ is $\Delta_0(M/L)$ (this is justified 
by \eqref{tors}). We put:
$$\vol^{(2)}(M)=\frac{\vol^{(2)}(L)}{\mathcal{M}(\Delta_0(M/L))}.$$
First of all, we have of course to check this is well-defined:

\begin{lem}
The real number $\frac{\vol^{(2)}(L)}{\mathcal{M}(\Delta_0(M/L))}$ does not 
depend on the free submodule of maximal rank $L\subset M$.
\end{lem}

\begin{proof}
Let $L_1,L_2$ be two free submodules of maximal rank in $M$. Then 
$L_1\cap L_2$ is a submodule of maximal rank, and thus contains a free 
submodule of maximal rank. Thus it is enough to prove that the $\ell^2$-volume 
is the same when computed for two free submodules of maximal rank 
$L\supset L'$. Let $Q$ be the coordinate matrix of some basis of $L'$ in a 
basis of $L$. We have:
$$\vol^{(2)}(L')=\det{}_{\mathcal{N}(G)}(Q)\vol^{(2)}(L)=\mathcal{M}(\Delta_0(L/L'))\vol^{(2)}(L)$$
On the other hand, we have a short exact sequence  
$0\rightarrow L/L'\rightarrow M/L'\rightarrow M/L\rightarrow 0$ which gives:
$$\mathcal{M}(\Delta_0(M/L'))=\mathcal{M}(\Delta_0(L/L'))\mathcal{M}(\Delta_0(M/L))$$
(by multiplicativity of $\Delta_0$ and the Mahler measure). Thus:
\begin{eqnarray*}
\frac{\vol^{(2)}(L')}{\mathcal{M}(\Delta_0(M/L'))}&=&\frac{\vol^{(2)}(L)}{\mathcal{M}(\Delta_0(M/L'))\mathcal{M}(\Delta_0(L/L'))}\\
       &=&\frac{\vol^{(2)}(L)}{\mathcal{M}(\Delta_0(M/L))}.
\end{eqnarray*}
This concludes the proof.
\end{proof}


\subsection{Metric rank formula}

We can now state the $\ell^2$ metric rank formula; the proof is similar to 
that of the finite-dimensional case.

\begin{lem}
Let $M,M'$ be $\mathbb{Z}[G]$-submodules of free modules and $f$ a 
morphism from $M$ to $M'$; we have then that:
\begin{equation}
\vol^{(2)}(f(M))=\frac{\det{}_{\mathcal{N}(G)}(f)\vol^{(2)}(M)}{\vol^{(2)}(\ker(f))}.
\label{rangmetr2}
\end{equation}
\label{rangmetrl}
\end{lem}

\begin{proof}
Renormalizing by an homothety we may suppose $\vol^{(2)}(M)=1$. 
Let $L'\subset\ker(f)$ and $L''\subset\ker(f)^{\bot}$ be free submodules 
of maximal rank; the direct sum $L'\oplus L''$ is then a free submodule of 
maximal rank in $M$. Put $T'=\ker(f)/L'$, $T''=\im(f)/f(L'')\cong 
M/\ker(f)\oplus L''$ and $T=M/L'\oplus L''$. We get a short exact sequence 
$0\rightarrow T'\rightarrow T\rightarrow T''\rightarrow0$, for example by 
applying the Snake Lemma to the diagram:
\[\begin{CD}
0@>>> \ker(f) @>>>     M      @>>> \im(f) @>>>0\\
&  &    @AAA         @AAA          @AA{f}A\\
0@>>>    L'    @>>> L'\oplus L'' @>>>   L''   @>>>0         
\end{CD}\]
It follows that:
$$\mathcal{M}(\Delta_0(T))=\mathcal{M}(\Delta_0(T''))\mathcal{M}(\Delta_0(T')).$$
On the other hand, since $L'$ and $L''$ are orthogonal we have:
$$\vol^{(2)}(L')\vol^{(2)}(L'')=\vol^{(2)}(L'+L'');$$
and since $M$ has unit volume this yields:
\begin{eqnarray*}
\vol^{(2)}(L')\vol^{(2)}(L'')&=&\mathcal{M}(\Delta_0(T))\\
                           &=&\mathcal{M}(\Delta_0(T''))\mathcal{M}(\Delta_0(T')).
\end{eqnarray*}
As $f_{|L''}$ is injective, $f(L'')$ is a free submodule of maximal rank in 
$\im(f)$ and thus:
\begin{eqnarray*}
\vol^{(2)}(\im(f))&=&\vol^{(2)}(f(L''))/\mathcal{M}(\Delta_0(\im(f)/f(L'')));\\
                  &=&\vol^{(2)}(f(L''))/\mathcal{M}(\Delta_0(T'')).
\end{eqnarray*}
Finally we have:
$$\vol^{(2)}(f(L''))=\det{}_{\mathcal{N}(G)}(f)\vol^{(2)}(L'').$$
From the three equalities above it follows that:
\begin{eqnarray*}
\vol^{(2)}(\im(f))&=&\frac{\det{}_{\mathcal{N}(G)}(f)\vol^{(2)}(L'')}{\mathcal{M}(\Delta_0(T''))}\\
                  &=&\frac{\det{}_{\mathcal{N}(G)}(f)\mathcal{M}(\Delta_0(T'))}{\vol^{(2)}(L')}\\
                  &=&\frac{\det{}_{\mathcal{N}(G)}(f)}{\vol^{(2)}(\ker(f))},
\end{eqnarray*}
and this concludes the proof of \eqref{rangmetr2}.
\end{proof}


\subsection{Computation of the $\ell^2$-torsion}

Let $C_*,d_*$ be a finite complex of free finitely generated 
$\mathbb{Z}[G]$-modules. 
Let $L$ be a free maximal rank submodule in $H_i(C)$; then $L$ lifts 
to a free submodule $L'$ in the orthogonal of $\im(d_{i+1})$ in $\ker(d_i)$. 
Let $T$ be the quotient $(H_i(C)/H_i(C)_{tors})/L$; we define:
$$\vol^{(2)}(H_i(C))=\mathcal{M}(\Delta_0(T))\vol^{(2)}(L').$$
In particular, when $C_*$ is $\ell^2$-acyclic this is equal to one 
since the orthogonal of $\im(d_{i+1})$ in $\ker(d_i)$ is zero.

\begin{lem}
With notation as above, the following equality holds:
\begin{equation*}
\tau^{(2)}(C)=\prod_i\left(\frac{\mathcal{M}(\Delta_{b_i^{(2)}(C)}(H_i(C))}{\vol^{(2)}(H_i(C))}\right)^{(-1)^i}
\end{equation*}
\label{calctors}
\end{lem}

\begin{proof}
We denote by $M_i$ the orthogonal of $\im(d_{i+1})$ in $\ker(d_i)$; this is 
a submodule of $C_i$ of rank $b_i^{(2)}(C)$. 
Let $L'$ be a maximal rank free submodule in $\im(d_{i-1})$ and $L''$ in 
$M_i$; put $L=L'\oplus L''$, which is a free maximal rank submodule 
in $\ker(d_i)$. Since  $M_i$ and $\im(d_{i-1})$ are orthogonal we 
get $\vol^{(2)}(L)=\vol^{(2)}(L')\vol^{(2)}(L'')$. Putting $T=\ker(d_i)/L$ 
it follows that:
$$\vol^{(2)}(\ker(d_i))=\frac{\vol^{(2)}(L)}{\mathcal{M}(\Delta_0(T))}=\frac{\vol^{(2)}(L')\vol^{(2)}(L'')}{\mathcal{M}(\Delta_0(T))}.$$
Put $T''=M_i/L''$ and $T'=\im(d_{i-1})/L'$; there is a natural 
injection of $T'\oplus T''$ into $T$, whose cokernel is naturally 
isomorphic to $H_i(C)/M_i=:T_0$. Thus we get an exact sequence  
$0\rightarrow T'\oplus T''\rightarrow T\rightarrow T_0 \rightarrow 0$, and 
the equality:
$$\mathcal{M}(\Delta_0(T))=\mathcal{M}(\Delta_0(T'))\mathcal{M}(\Delta_0(T''))\mathcal{M}(\Delta_0(T_0)).$$
Finally, putting $T_1=\widehat{H}_i(C)/M_i$ we get the exact 
sequence 
$0\rightarrow T_1\rightarrow T_0\rightarrow H_i(C)_{tors}\rightarrow0$ and we 
can compute:
$$\mathcal{M}(\Delta_0(T_0))=\mathcal{M}(\Delta_0(T_1))\mathcal{M}(\Delta_0(H_i(C)_{tors})).$$
Putting the three equalities above together gives:
\begin{eqnarray*}
\vol^{(2)}(\ker(d_i))&=&\frac{\vol^{(2)}(L')}{\mathcal{M}(\Delta_0(T'))}\times\frac{\vol^{(2)}(L'')}{\mathcal{M}(\Delta_0(T''))}\times\frac{1}{\mathcal{M}(\Delta_0(T_0))}\\
                     &=&\vol^{(2)}(\im(d_{i-1}))\times\frac{\vol^{(2)}(M_i)}{\mathcal{M}(\Delta_0(T_1))}\times\frac{1}{\mathcal{M}(\Delta_0(H_i(C)_{tors}))}\\
                     &=&\frac{\vol^{(2)}(\im(d_{i-1}))\vol^{(2)}(\widehat{H}_i(C))}{\mathcal{M}(\Delta_0(H_i(C)_{tors}))}.
\end{eqnarray*}
From this and \eqref{rangmetr2} it follows that:
$$\vol^{(2)}(\ker(d_i))=\frac{\det{}_{\mathcal{N}(G)}(d_{i-1})\vol^{(2)}(\widehat{H}_i(C))}{\mathcal{M}(\Delta_0(H_i(C)_{tors}))\vol^{(2)}(\ker d_{i-1})}.$$
Taking the alternating product over $i$ we obtain:
$$1=\prod_i \left(\frac{\det{}_{\mathcal{N}(G)}(d_{i-1})\vol^{(2)}(\widehat{H}_i(C))}{\mathcal{M}(\Delta_0(H_i(C)_{tors}))}\right)^{(-1)^i}.$$
\end{proof}


\bibliographystyle{amsalpha}
\bibliography{biblio}

\providecommand{\bysame}{\leavevmode\hbox to3em{\hrulefill}\thinspace}
\providecommand{\MR}{\relax\ifhmode\unskip\space\fi MR }
\providecommand{\MRhref}[2]{%
  \href{http://www.ams.org/mathscinet-getitem?mr=#1}{#2}
}
\providecommand{\href}[2]{#2}
\begin{thebibliography}{CFK11}

\bibitem[BV]{Bergeron/Venkatesh}
Nicolas Bergeron and Akshay Venkatesh, \emph{The asymptotic growth of torsion
  homology for arithmetic groups}, To appear in JIMJ, math.NT/1004.1083v1.

\bibitem[CFK11]{KFC}
Jae~Choon Cha, Stefan Friedl, and Taehee Kim, \emph{The cobordism group of
  homology cylinders}, Compos. Math. \textbf{147} (2011), no.~3, 914--942.
  \MR{2801405}

\bibitem[CW03]{Clair/Whyte}
Bryan Clair and Kevin Whyte, \emph{Growth of {B}etti numbers}, Topology
  \textbf{42} (2003), no.~5, 1125--1142. \MR{MR1978050 (2004e:58057)}

\bibitem[EW99]{Everest/Ward}
Graham Everest and Thomas Ward, \emph{Heights of polynomials and entropy in
  algebraic dynamics}, Universitext, Springer-Verlag London Ltd., London, 1999.
  \MR{MR1700272 (2000e:11087)}

\bibitem[GAS91]{G-A/S}
Francisco Gonz{\'a}lez-Acu{\~n}a and Hamish Short, \emph{Cyclic branched
  coverings of knots and homology spheres}, Rev. Mat. Univ. Complut. Madrid
  \textbf{4} (1991), no.~1, 97--120. \MR{MR1142552 (93g:57004)}

\bibitem[Gor72]{Gordon}
C.~McA. Gordon, \emph{Knots whose branched cyclic coverings have periodic
  homology}, Trans. Amer. Math. Soc. \textbf{168} (1972), 357--370.
  \MR{MR0295327 (45 \#4394)}

\bibitem[Hat02]{Hatcher}
Allen Hatcher, \emph{Algebraic topology}, Cambridge University Press,
  Cambridge, 2002. \MR{1867354 (2002k:55001)}

\bibitem[Hem76]{Hempel}
John Hempel, \emph{{$3$}-{M}anifolds}, Princeton University Press, Princeton,
  N. J., 1976, Ann. of Math. Studies, No. 86. \MR{MR0415619 (54 \#3702)}

\bibitem[Kow08]{Kowalski}
Emmanuel Kowalski, \emph{The large sieve and its applications}, Cambridge
  Tracts in Mathematics, vol. 175, Cambridge University Press, Cambridge, 2008,
  Arithmetic geometry, random walks and discrete groups. \MR{2426239
  (2009f:11123)}

\bibitem[Lac10]{Lackenby}
Marc Lackenby, \emph{Finite covering spaces of 3-manifolds}, Proceedings of the
  {I}nternational {C}ongress of {M}athematicians, {H}yberabad, {I}ndia, 2010.

\bibitem[Lam06]{Lam}
T.~Y. Lam, \emph{Serre's problem on projective modules}, Springer Monographs in
  Mathematics, Springer-Verlag, Berlin, 2006. \MR{MR2235330 (2007b:13014)}

\bibitem[Law83]{Lawton}
Wayne~M. Lawton, \emph{A problem of {B}oyd concerning geometric means of
  polynomials}, J. Number Theory \textbf{16} (1983), no.~3, 356--362.
  \MR{MR707608 (84i:10056)}

\bibitem[Le]{Le}
Thang Le, \emph{Homology torsion growth and {M}ahler measure},
  math.GT/1010.4199.

\bibitem[Lic97]{Lickorish}
W.~B.~Raymond Lickorish, \emph{An introduction to knot theory}, Graduate Texts
  in Mathematics, vol. 175, Springer-Verlag, New York, 1997. \MR{MR1472978
  (98f:57015)}

\bibitem[LLS11]{LLS}
Peter Linnell, Wolfgang L{\"u}ck, and Roman Sauer, \emph{The limit of {$\Bbb
  F_p$}-{B}etti numbers of a tower of finite covers with amenable fundamental
  groups}, Proc. Amer. Math. Soc. \textbf{139} (2011), no.~2, 421--434.
  \MR{2736326 (2012a:55016)}

\bibitem[L{\"u}c94]{Luck2}
W.~L{\"u}ck, \emph{Approximating {$L^2$}-invariants by their finite-dimensional
  analogues}, Geom. Funct. Anal. \textbf{4} (1994), no.~4, 455--481.
  \MR{MR1280122 (95g:58234)}

\bibitem[L{\"u}c02]{Luck}
Wolfgang L{\"u}ck, \emph{{$L^2$}-invariants: theory and applications to
  geometry and {$K$}-theory}, Ergebnisse der Mathematik und ihrer Grenzgebiete.
  3. Folge., vol.~44, Springer-Verlag, Berlin, 2002. \MR{MR1926649
  (2003m:58033)}

\bibitem[MM82]{Mayberry/Murasugi}
John~P. Mayberry and Kunio Murasugi, \emph{Torsion-groups of abelian coverings
  of links}, Trans. Amer. Math. Soc. \textbf{271} (1982), no.~1, 143--173.
  \MR{MR648083 (84d:57004)}

\bibitem[Por04]{Porti}
Joan Porti, \emph{Mayberry-{M}urasugi's formula for links in homology
  3-spheres}, Proc. Amer. Math. Soc. \textbf{132} (2004), no.~11, 3423--3431
  (electronic). \MR{MR2073320 (2005e:57004)}

\bibitem[Ril90]{Riley}
Robert Riley, \emph{Growth of order of homology of cyclic branched covers of
  knots}, Bull. London Math. Soc. \textbf{22} (1990), no.~3, 287--297.
  \MR{MR1041145 (92g:57017)}

\bibitem[Sch95]{Schmidt}
Klaus Schmidt, \emph{Dynamical systems of algebraic origin}, Progress in
  Mathematics, vol. 128, Birkh\"auser Verlag, Basel, 1995. \MR{MR1345152
  (97c:28041)}

\bibitem[Sch00]{Schinzel}
A.~Schinzel, \emph{Polynomials with special regard to reducibility},
  Encyclopedia of Mathematics and its Applications, vol.~77, Cambridge
  University Press, Cambridge, 2000, With an appendix by Umberto Zannier.
  \MR{1770638 (2001h:11135)}

\bibitem[Ser00]{Serre}
Jean-Pierre Serre, \emph{Local algebra}, Springer Monographs in Mathematics,
  Springer-Verlag, Berlin, 2000, Translated from the French by CheeWhye Chin
  and revised by the author. \MR{MR1771925 (2001b:13001)}

\bibitem[SW02a]{Silver/Williams1}
Daniel~S. Silver and Susan~G. Williams, \emph{Mahler measure, links and
  homology growth}, Topology \textbf{41} (2002), no.~5, 979--991. \MR{MR1923995
  (2003h:57011)}

\bibitem[SW02b]{Silver/Williams2}
\bysame, \emph{Torsion numbers of augmented groups with applications to knots
  and links}, Enseign. Math. (2) \textbf{48} (2002), no.~3-4, 317--343.
  \MR{MR1955605 (2003m:57025)}

\bibitem[SW09]{Silver/Williams3}
\bysame, \emph{Dynamics of twisted {A}lexander invariants}, Topology Appl.
  \textbf{156} (2009), no.~17, 2795--2811. \MR{2556037}

\bibitem[Tur01]{Turaev}
Vladimir Turaev, \emph{Introduction to combinatorial torsions}, Lectures in
  Mathematics ETH Z\"urich, Birkh\"auser Verlag, Basel, 2001, Notes taken by
  Felix Schlenk. \MR{MR1809561 (2001m:57042)}

\end{thebibliography}

\end{document}